\documentclass[12pt,a4paper,draft]{amsart}

\usepackage{hyperref} 
\usepackage[all]{xy}

\usepackage{amssymb,amsrefs}

\newcommand{\bbA}{\mathbb{A}} %
\newcommand{\bbD}{\mathbb{D}} %
\newcommand{\bbE}{\mathbb{E}} %
\newcommand{\bbM}{\mathbb{M}} %
\newcommand{\bbN}{\mathbb{N}} %
\newcommand{\bbX}{\mathbb{X}} %
\newcommand{\bbZ}{\mathbb{Z}} %

\newcommand{\calA}{\mathcal{A}} %
\newcommand{\calC}{\mathcal{C}} %
\newcommand{\calH}{\mathcal{H}} %
\newcommand{\calI}{\mathcal{I}} %
\newcommand{\calO}{\mathcal{O}} %
\newcommand{\calP}{\mathcal{P}} %
\newcommand{\calQ}{\mathcal{Q}} %
\newcommand{\calR}{\mathcal{R}} %
\newcommand{\calU}{\mathcal{U}} %
\newcommand{\calV}{\mathcal{V}} %
\newcommand{\calX}{\mathcal{X}} %
\newcommand{\calY}{\mathcal{Y}} %

\newcommand{\bd}{\mathbf{d}}
\newcommand{\be}{\mathbf{e}}
\newcommand{\bh}{\mathbf{h}}
\newcommand{\bx}{\mathbf{x}}
\newcommand{\bP}{\mathbf{P}}
\newcommand{\bQ}{\mathbf{Q}}
\newcommand{\bR}{\mathbf{R}}
\newcommand{\bone}{\mathbf{1}}
\newcommand{\bDelta}{\mathbf{\Delta}}

\let\mod=\undefined
\DeclareMathOperator{\GL}{GL} %
\DeclareMathOperator{\SI}{SI} %
\DeclareMathOperator{\add}{add} %
\DeclareMathOperator{\Ext}{Ext} %
\DeclareMathOperator{\Hom}{Hom} %
\DeclareMathOperator{\ind}{ind} %
\DeclareMathOperator{\Ker}{Ker} %
\DeclareMathOperator{\mod}{mod} %
\DeclareMathOperator{\rep}{rep} %
\DeclareMathOperator{\idim}{idim} %
\DeclareMathOperator{\pdim}{pdim} %
\DeclareMathOperator{\codim}{codim} %
\DeclareMathOperator{\Coker}{Coker} %
\DeclareMathOperator{\gldim}{gl.dim} %
\DeclareMathOperator{\bdim}{\mathbf{dim}} %

\newcommand{\ol}{\overline}

\newcounter{mainclaim}
\newtheorem{maintheorem}[mainclaim]{Theorem}
\newtheorem{maincorollary}[mainclaim]{Corollary}

\newcounter{claim}[section]

\newtheorem{corollary}[claim]{Corollary}
\newtheorem{lemma}[claim]{Lemma}
\newtheorem{proposition}[claim]{Proposition}

\title[Maximal orbits for Euclidean quivers]%
{Normality of maximal orbit closures \\ for Euclidean quivers}

\author{Grzegorz Bobi\'nski}

\address{Faculty of Mathematics and Computer Science \\ Nicolaus
Copernicus University \\ ul.~Chopina 12/18 \\ 87-100 Toru\'n \\
Poland}

\email{gregbob@mat.umk.pl}


\begin{document}

\begin{abstract}
Let $\Delta$ be an Euclidean quiver. We prove that the closures of
the maximal orbits in the varieties of representations of $\Delta$
are normal and Cohen--Macaulay (even complete intersections).
Moreover, we give a generalization of this result for the tame
concealed-canonical algebras.
\end{abstract}

\maketitle

\section*{Introduction and the main results}

Throughout the paper $k$ is a fixed algebraically closed field. By
$\bbZ$, $\bbN$ and $\bbN_+$ we denote the sets of the integers,
the non-negative integers and the positive integers, respectively.
Finally, if $i, j \in \bbZ$, then $[i, j] := \{ l \in \bbZ \mid i
\leq l \leq j \}$ (in particular, $[i, j] = \varnothing$ if $i >
j$).

Let $A$ be a finite dimensional $k$-algebra. Given a non-negative
integer $d$ one defines $\mod_A (d)$ as the set of all $k$-algebra
homomorphisms from $A$ to the algebra $\bbM_{d \times d} (k)$ of
$d \times d$-matrices. This set has a structure of an affine
variety and its points represent $d$-dimensional $A$-modules.
Consequently, we call $\mod_A (d)$ the variety of $A$-modules of
dimension $d$. The general linear group $\GL (d)$ acts on $\mod_A
(d)$ by conjugation: $(g \cdot m) (a) := g m (a) g^{-1}$ for $g
\in \GL (d)$, $m \in \mod_A (d)$ and $a \in A$. The orbits with
respect to this action are in one-to-one correspondence with the
isomorphism classes of the $d$-dimensional $A$-modules. Given a
$d$-dimensional $A$-module $M$ we denote the orbit in $\mod_A (d)$
corresponding to the isomorphism class of $M$ by $\calO (M)$ and
its Zariski-closure by $\ol{\calO (M)}$.

Singularities appearing in the orbit closures $\ol{\calO (M)}$ for
modules $M$ over an algebra $A$ are an object of intensive studies
(see for example ~\cites{AbeasisDelFraKraft1981,
BenderBongartz2003, BobinskiZwara2006, Bongartz1994, LocZwara,
SkowronskiZwara2003, Zwara2005a, Zwara2005b, Zwara2006,
Zwara2007}, we also refer to a survey article of
Zwara~\cite{Zwara2011}). In particular, Zwara and the
author~\cite{BobinskiZwara2002} proved that if $A$ is a hereditary
algebra of Dynkin type $\bbA$ or $\bbD$, then $\ol{\calO (M)}$ is
a normal Cohen--Macaulay variety, which has rational singularities
if the characteristic of $k$ is $0$. Recall, that
Gabriel~\cite{Gabriel1972} proved that the hereditary algebras of
Dynkin type are precisely the hereditary algebras of finite
representation type. Thus, it an interesting question if the orbit
closures have good geometric properties for all hereditary
algebras of finite representation type. The remaining case of
hereditary algebras of type $\bbE$ is still open, but there are
some partial results in this direction~\cite{Zwara2002}. On the
other hand, Zwara~\cite{Zwara2003} exhibited an example of a
module over the Kronecker algebra whose orbit closure is neither
normal nor Cohen--Macaulay . This example generalizes easily to an
arbitrary hereditary algebra of infinite representation
type~\cite{Chindris2007}. However, it is still an interesting
problem to determine for which classes of modules over hereditary
algebras of infinite representation type, the corresponding orbit
closures have good properties. In the paper, we study modules $M$
such that $\calO (M)$ is maximal, i.e.\ there is no module $N$
such that $\calO (M) \subseteq \ol{\calO (N)}$ and $\calO (M) \neq
\calO (N)$.

According to famous Drozd's Tame and Wild
Theorem~\cites{Drozd1980, CrawleyBoevey1988} the finite
dimensional algebras of infinite representation type can be
divided into two disjoint classes. One class consists of the tame
algebras, for which the indecomposable modules of a given
dimension form a finite number of one-parameter families. The
other class consists of the wild algebras, for which the
classification of the indecomposable modules is as complicated as
the classification of two non-commuting endomorphisms of a finite
dimensional vector space, hence is considered to be hopeless.
There are examples showing that varieties of modules over tame
algebras have often better properties than those over wild
algebras (see for example~\cites{BobinskiSkowronski1999,
Chindris2009, SkowronskiWeyman2000, SkowronskiZwara1998}).
Consequently, we concentrate in the paper on the maximal orbits
over the tame hereditary algebras. We recall that the tame
hereditary algebras are precisely the hereditary algebras of
Euclidean type.

The following theorem is the main result of the paper.

\begin{maintheorem} \label{main theorem}
Let $M$ be a module over a tame hereditary algebra. If $\calO (M)$
is maximal, then $\ol{\calO (M)}$ is a normal complete
intersection \textup{(}in particular, Cohen--Macaulay\textup{)}.
\end{maintheorem}

It is known (see for example~\cite{Ringel1980}*{Corollary~3.6})
that $\calO (M)$ is maximal for each indecomposable module over a
tame hereditary algebra. Consequently, we get the following.

\begin{maincorollary}
If $M$ is an indecomposable module over a tame hereditary algebra,
then $\ol{\calO (M)}$ is a normal complete intersection
\textup{(}in particular, Cohen--Macaulay\textup{)}.
\end{maincorollary}

Now we present the strategy of the proof of Theorem~\ref{main
theorem}. Let $M$ be a module over a tame hereditary algebra $A$
such that $\calO (M)$ is maximal. If $\Ext_A^1 (M, M) = 0$, then
it is well known that $\ol{\calO (M)}$ is smoothly equivalent to
an affine space, hence the claim is obvious in this case. Thus we
may concentrate on the case $\Ext_A^1 (M, M) \neq 0$. It follows
from~\cite{Ringel1980}*{proof of Corollary~3.6} that in this
situation $M$ is periodic with respect to the action of the
Auslander--Reiten translation $\tau$. Consequently,
Theorem~\ref{main theorem} follows from the following.

\begin{maintheorem} \label{main theorem bis}
Let $M$ be a $\tau$-periodic module over a tame hereditary
algebra. If $\calO (M)$ is maximal, then $\ol{\calO (M)}$ is a
complete intersection \textup{(}in particular,
Cohen--Macaulay\textup{)}.
\end{maintheorem}

If $A$ is a tame hereditary algebra, then the $\tau$-periodic
$A$-modules are direct sums of indecomposable modules, which lie
in the sincere separating family of tubes in the Auslander--Reiten
quiver of $A$. Existence of such families charecterizes the
concealed-canonical algebras~\cites{LenzingdelaPena1999,
Skowronski1996}. Recall~\cite{LenzingMeltzer1996} that an algebra
$A$ is called concealed-canonical if there exists a tilting bundle
over a weighted projective line whose endomorphism ring is
isomorphic to $A$. Thus it is natural to try to generalize
Theorem~\ref{main theorem bis} to the case of tame
concealed-canonical algebras. Before we formulate this
generalization, we present necessary definitions.

Let $A$ be a tame concealed-canonical algebra. For an $A$-module
$M$ we denote by $\bdim M$ its dimension vector, i.e.\ the
sequence indexed by the isomorphism classes of the simple
$A$-modules, which counts the multiplicities of the composition
factors in the Jordan--H\"older filtration of $M$. In general, a
sequence of non-negative integers indexed by the isomorphism
classes of the simple $A$-modules is called a dimension vector. We
call a dimension vector $\bd$ singular if $\langle \bd, \bd
\rangle_A = 0$ and there exists a dimension vector $\bx$ such that
$\bx \leq \bd$, $\langle \bx, \bx \rangle_A = 0$ and $|\langle
\bx, \bd \rangle_A| = 2$, where $\langle -, - \rangle_A$ denotes
the corresponding homological bilinear form (see
Section~\ref{section quivers}). In Proposition~\ref{proposition
singular} we describe the tame concealed-canonical algebras for
which there exist singular dimension vectors. In particular, this
description implies that singular dimension vectors do not exist
for the tame hereditary algebras.

We have the following generalization of Theorem~\ref{main theorem
bis}.

\begin{maintheorem} \label{main theorem prim}
Let $M$ be a $\tau$-periodic module over a tame
concealed-canonical algebra such that $\calO (M)$ is maximal. Then
$\ol{\calO (M)}$ is a complete intersection \textup{(}in
particular, Cohen--Macaulay\textup{)}. Moreover, $\ol{\calO (M)}$
is not normal if and only if $\bdim M$ is singular and $\tau M
\simeq M$.
\end{maintheorem}

In the paper we concentrate on the proof of Theorem~\ref{main
theorem prim}. Instead of using the framework of modules over
algebras and the corresponding varieties, we use the framework of
representations of quivers (and the corresponding varieties).
Gabriel's Theorem~\cite{Gabriel1972} says that we may do this
replacement on the level of modules and representations, while a
result of Bongartz~\cite{Bongartz1991} justifies this passage on
the level of varieties. For the background on the representation
theory we refer to~\cites{AssemSimsonSkowronski2006,
SimsonSkowronski2007a, SimsonSkowronski2007b}.

The paper is organized as follows. In Section~\ref{section
quivers} we recall basic information about quivers and their
representations. Next, in Section~\ref{section tubular} we gather
facts about the categories of modules over the tame
concealed-canonical algebras. In Section~\ref{section varieties}
we introduce varieties of representations of quivers, while in
Section~\ref{section semiinvariants} we review facts on
semi-invariants with particular emphasis on the case of tame
concealed-canonical algebras. Next, in Section~\ref{section
lemmas} we present a series of facts, which we later use in
Sections~\ref{section nonsingular} and~\ref{section singular} to
study orbit closures for the non-singular and singular dimension
vectors, respectively. Moreover, in Section~\ref{section singular}
we make a remark about relationship between the degenerations and
the hom-order for the tame concealed-canonical algebras. Finally,
in Section~\ref{section proof} we give the proof of
Theorem~\ref{main theorem prim}.

The author acknowledges the support from the Research Grant No.\ N
N201 269135 of the Polish Ministry of Science and Higher
Education.

\section{Quivers and their representations} \label{section
quivers}

By a quiver $\Delta$ we mean a finite set $\Delta_0$ (called the
set of vertices of $\Delta$) together with a finite set $\Delta_1$
(called the set of arrows of $\Delta$) and two maps $s, t :
\Delta_1 \to \Delta_0$, which assign to each arrow $\alpha$ its
starting vertex $s \alpha$ and terminating vertex $t \alpha$,
respectively. By a path of length $n \in \bbN_+$ in a quiver
$\Delta$ we mean a sequence $\sigma = (\alpha_1, \ldots,
\alpha_n)$ of arrows such that $s \alpha_i = t \alpha_{i + 1}$ for
each $i \in [1, n - 1]$. In particular, we treat every arrow in
$\Delta$ as a path of length $1$. In the above situation we put
$\ell \sigma := n$, $s \sigma := s \alpha_n$ and $t \sigma := t
\alpha_1$. Moreover, for each vertex $x$ we have a trivial path
$\bone_x$ at $x$ such that $\ell \bone_x := 0$ and $s \bone_x := x
=: t \bone_x$. A subquiver $\Delta'$ of a quiver $\Delta$ is
called convex if $\alpha_i \in \Delta_1'$ for each $i \in [1, n]$,
provided $(\alpha_1, \ldots, \alpha_n)$ is a path in $\Delta$ such
that $t \alpha_1, s \alpha_n \in \Delta_0'$.

For the rest of the paper we assume that the considered quivers do
not have oriented cycles, where by an oriented cycle we mean a
path $\sigma$ of positive length such that $s \sigma = t \sigma$.

Let $\Delta$ be a quiver. We define its path category $k \Delta$
to be the category whose objects are the vertices of $\Delta$ and,
for $x, y \in \Delta_0$, the morphisms from $x$ to $y$ are the
formal $k$-linear combinations of paths starting at $x$ and
terminating at $y$. For $x, y \in \Delta_0$ we denote by $k \Delta
(x, y)$ the space of the morphisms from $x$ to $y$ in $k \Delta$.
If $\omega \in k \Delta (x, y)$ for $x, y \in \Delta_0$, then we
write $s \omega := x$ and $t \omega := y$. By a representation of
$\Delta$ we mean a functor from $k \Delta$ to the category $\mod
k$ of finite dimensional vector spaces. We denote the category of
the representations of $\Delta$ by $\rep \Delta$. Observe that
every representation of $\Delta$ is uniquely determined by its
values on the vertices and the arrows. Given a representation $M$
of $\Delta$ we denote by $\bdim M$ its dimension vector defined by
$(\bdim M) (x) := \dim_k M (x)$ for $x \in \Delta_0$. Observe the
$\bdim M \in \bbN^{\Delta_0}$ for each representation $M$ of
$\Delta$. We call the elements of $\bbN^{\Delta_0}$ dimension
vectors. A dimension vector $\bd$ is called sincere if $\bd (x)
\neq 0$ for each $x \in \Delta_0$.

By a relation in a quiver $\Delta$ we mean a $k$-linear
combination of paths of lengths at least $2$ having a common
starting vertex and a common terminating vertex. Note that each
relation in a quiver $\Delta$ is a morphism in $k \Delta$. A set
$R$ of relations in a quiver $\Delta$ is called minimal if
$\langle R \setminus \{ \rho \} \rangle \neq \langle R \rangle$
for each $\rho \in R$, where for a set $X$ of morphisms in
$\Delta$ we denote by $\langle X \rangle$ the ideal in $k \Delta$
generated by $X$. Observe that each minimal set of relations is
finite. By a bound quiver $\bDelta$ we mean a quiver $\Delta$
together with a minimal set $R$ of relations. Given a bound quiver
$\bDelta$ we denote by $k \bDelta$ its path category, i.e.\ $k
\bDelta := k \Delta / \langle R \rangle$. Moreover, for $x, y \in
\Delta_0$ we denote by $k \bDelta (x, y)$ the space of the
morphisms from $x$ to $y$ in $k \bDelta$. By a representation of a
bound quiver $\bDelta$ we mean a functor from $k \bDelta$ to $\mod
k$. In other words, a representation of $\bDelta$ is a
representation $M$ of $\Delta$ such that $M (\rho) = 0$ for each
$\rho \in R$. We denote the category of the representations of a
bound quiver $\bDelta$ by $\rep \bDelta$. Moreover, we denote by
$\ind \bDelta$ the full subcategory of $\rep \bDelta$ consisting
of the indecomposable representations. It is known that $\rep
\bDelta$ is an abelian Krull--Schmidt category. A bound quiver
$\bDelta'$ is called a convex subquiver of a bound quiver
$\bDelta$ if $\Delta'$ is a convex subquiver of $\Delta$ and $R' =
R \cap k \Delta'$. If $\bDelta'$ is a convex subquiver of a bound
quiver $\bDelta$, then $\rep \bDelta'$ can be naturally identified
with an exact subcategory of $\rep \bDelta$, where by an exact
subcategory of $\rep \bDelta$ we mean a full subcategory $\calX$
of $\rep \bDelta$ such that $\calX$ is an abelian category and the
inclusion functor $\calX \hookrightarrow \rep \bDelta$ is exact.
In particular, if $\bDelta'$ is a convex subcategory of a tame
bound quiver $\bDelta$, then $\bDelta'$ is either tame or
representation-finite (we say that a bound quiver $\bDelta$ is
tame/representation-finite if $\rep \bDelta$ is of tame/finite
representation type, respectively).

Let $\bDelta$ be a bound quiver. For each vertex $x$ of $\Delta$
we denote by $S_x$ the simple representation at $x$, i.e.\ $S_x
(x) := k$, $S_x (y) := 0$ for $y \in \Delta_0 \setminus \{ x \}$,
and $S_x (\alpha) := 0$ for $\alpha \in \Delta_1$. More generally,
if $\bd$ is a dimension vector, then we put $S^{\bd} :=
\bigoplus_{x \in \Delta_0} S_x^{\bd (x)}$. Next, for each vertex
$x$ we denote by $P_x$ the projective representation at $x$
defined in the following way: $P_x (y) := k \bDelta (x, y)$ for $y
\in \Delta_0$ and $P_x (\omega)$ is the composition (on the left)
with $\omega$ for a morphism $\omega$ in $k \bDelta$. If $M$ is a
representation of $\bDelta$ and $x \in \Delta_0$, then the map
\[
\Hom_{\bDelta} (P_x, M) \to M (x), \; f \mapsto f (\bone_x),
\]
is an isomorphism. In particular, this implies that
\[
\Hom_{\bDelta} (P_x, P_y) \simeq k \bDelta (y, x)
\]
for any $x, y \in \Delta_0$. If $\omega \in k \bDelta (y, x)$, we
denote the corresponding map $P_x \to P_y$ by $P_\omega$. Observe
that $P_\omega$ is the composition (on the right) with $\omega$.
Moreover, if $M$ is a representation of $\bDelta$, then, under the
above isomorphisms, $\Hom_{\bDelta} (P_\omega, M)$ equals $M
(\omega)$.

Let $\bDelta$ be a bound quiver. If $P_1 \xrightarrow{f} P_0 \to M
\to 0$ is a (fixed) minimal projective presentation of a
representation $M$ of $\bDelta$, then we put
\[
\tau M := \Hom_k (\Coker \Hom_{\bDelta} (f, \bigoplus_{x \in
\Delta_0} P_x), k).
\]
We define $\tau^- M$ dually. Note that $\tau M = 0$ ($\tau^- M =
0$) if and only if $M$ is projective (injective, respectively).
Moreover, $\tau^- \tau X \simeq X$ ($\tau \tau^- X \simeq X$) for
each indecomposable representation $X$ of $\bDelta$, which is not
projective (injective, respectively). We say that a representation
$M$ of $\bDelta$ is periodic if there exists $n \in \bbN_+$ such
that $\tau^n M \simeq M$. We have a celebrated Auslander--Reiten
formula~\cite{AssemSimsonSkowronski2006}*{Section~IV.2}, which
implies that
\[
\dim_k \Ext_{\bDelta}^1 (M, N) = \dim_k \Hom_{\bDelta} (N, \tau M)
\]
for any representations $M$ and $N$ of $\bDelta$ such that
$\pdim_{\bDelta} M \leq 1$. Dually, if $M$ and $N$ are
representations of $\bDelta$ and $\idim_{\bDelta} N \leq 1$, then
\[
\dim_k \Ext_{\bDelta}^1 (M, N) = \dim_k \Hom_{\bDelta} (\tau^- N,
M).
\]

Let $\bDelta$ be a bound quiver. We define the corresponding Tits
forms $\langle -, - \rangle_{\bDelta} : \bbZ^{\Delta_0} \times
\bbZ^{\Delta_0} \to \bbZ$ and $q_{\bDelta} : \bbZ^{\Delta_0} \to
\bbZ$ by
\[
\langle \bd', \bd'' \rangle_{\bDelta} := \sum_{x \in \Delta_0}
\bd' (x) \bd'' (x) - \sum_{\alpha \in \Delta_1} \bd' (s \alpha)
\bd'' (t \alpha) + \sum_{\rho \in R} \bd' (s \rho) \bd'' (t \rho)
\]
for $\bd', \bd'' \in \bbZ^{\Delta_0}$, and $q_{\bDelta} (\bd) :=
\langle \bd, \bd \rangle_{\bDelta}$ for $\bd \in \bbZ^{\Delta_0}$.
Bongartz~\cite{Bongartz1983}*{Proposition~2.2} proved that
\begin{multline*}
\langle \bdim M, \bdim N \rangle_{\bDelta}
\\
= \dim_k \Hom_{\bDelta} (M, N) - \dim_k \Ext_{\bDelta}^1 (M, N) +
\dim_k \Ext_{\bDelta}^2 (M, N)
\end{multline*}
for any $M, N \in \rep \bDelta$, provided $\gldim \bDelta \leq 2$.

\section{Separating exact subcategories} \label{section tubular}

In this section we present facts about sincere separating exact
subcategories, which we use in our considerations. For the proofs
we refer to~\cites{Ringel1984, LenzingdelaPena1999}.

Let $\bDelta$ be a bound quiver and $\calX$ a full subcategory of
$\ind \bDelta$. We denote by $\add \calX$ the full subcategory of
$\rep \bDelta$ formed by the direct sums of representations from
$\calX$. We say that $\calX$ is an exact subcategory of $\ind
\bDelta$ if $\add \calX$ is an exact subcategory of $\rep
\bDelta$. We put
\begin{gather*}
\calX_+ := \{ X \in \ind \bDelta : \text{$\Hom_{\bDelta} (\calX,
X) = 0$} \}
\\
\intertext{and} %
\calX_- := \{ X \in \ind \bDelta : \text{$\Hom_{\bDelta} (X,
\calX) = 0$} \}.
\end{gather*}

Let $\bDelta$ be a bound quiver.
Following~\cite{LenzingdelaPena1999} we say that $\calR$ is a
sincere separating exact subcategory of $\ind \bDelta$ provided
the following conditions are satisfied:
\begin{enumerate}

\item
$\calR$ is an exact subcategory of $\ind \bDelta$ consisting of
periodic representations.

\item
$\ind \bDelta = \calR_+ \cup \calR \cup \calR_-$.

\item
$\Hom_{\bDelta} (X, \calR) \neq 0$ for each $X \in \calR_+$ and
$\Hom_{\bDelta} (\calR, X) \neq 0$ for each $X \in \calR_-$.

\item
$P \in \calR_+$ for each indecomposable projective representation
$P$ of $\bDelta$ and $I \in \calR_-$ for each indecomposable
injective representation $I$ of $\bDelta$.

\end{enumerate}
Lenzing and de la Pe\~na~\cite{LenzingdelaPena1999} proved that
there exists a sincere separating exact subcategory $\calR$ of
$\rep \bDelta$ if and only if $\bDelta$ is concealed-canonical,
i.e.\ $\rep \bDelta$ is equivalent to the category of modules over
a concealed-canonical algebra. In particular, if this the case,
then $\gldim \bDelta \leq 2$.

For the rest of the section we fix a bound quiver $\bDelta$ and a
sincere separating exact subcategory $\calR$ of $\ind \bDelta$.
Moreover, we put $\calP := \calR_+$ and $\calQ := \calR_-$.
Finally, we denote by $\bP$, $\bR$ and $\bQ$ the dimension vectors
of the representations from $\add \calP$, $\add \calR$ and $\add
\calQ$, respectively.

It is known that $\pdim_{\bDelta} P \leq 1$ for each $P \in \calP$
and $\idim_{\bDelta} Q \leq 1$ for each $Q \in \calQ$. Next,
$\pdim_{\bDelta} R = 1$ and $\idim_{\bDelta} R = 1$ for each $R
\in \calR$. Moreover, $\Hom_{\bDelta} (\calQ, \calP) = 0$. Since
the categories $\calP$ and $\calQ$ are closed under the actions of
$\tau$ and $\tau^-$, using the Auslander--Reiten formulas we also
obtain that $\Ext_{\bDelta}^1 (\calP, \calR \cup \calQ) = 0 =
\Ext_{\bDelta}^1 (\calP \cup \calR, \calQ)$. The above properties
imply that $\langle \bd', \bd'' \rangle_{\bDelta} \geq 0$ if
either $\bd' \in \bP$ and $\bd'' \in \bR + \bQ$ or $\bd' \in \bP +
\bR$ and $\bd'' \in \bQ$. Similarly, $\langle \bd'', \bd'
\rangle_{\bDelta} \leq 0$ if either $\bd' \in \bP$ and $\bd'' \in
\bR$ or $\bd' \in \bR$ and $\bd'' \in \bQ$.

We have $\calR = \coprod_{\lambda \in \bbX} \calR_\lambda$ for
some infinite set $\bbX$ and connected uniserial categories
$\calR_\lambda$, $\lambda \in \bbX$. For $\lambda \in \bbX$ we
denote by $r_\lambda$ the number of the pairwise non-isomorphic
simple objects in $\add \calR_\lambda$. Then $r_\lambda < \infty$.
Let $\bbX_0 := \{ \lambda \in \bbX : \text{$r_\lambda > 1$} \}$.
Then $|\bbX_0| < \infty$ and we call the sequence
$(r_\lambda)_{\lambda \in \bbX_0}$ the type of $\bDelta$ (this
definition does not depend on the choice of a sincere separating
exact subcategory of $\ind \bDelta$). It is known that $\bDelta$
is tame if and only if $\sum_{\lambda \in \bbX_0}
\frac{1}{r_\lambda} \geq |\bbX_0| - 2$, where by definition the
empty sum equals $0$. Observe that this implies that $|\bbX_0|
\leq 4$ provided $\bDelta$ is tame. Moreover, if $\bDelta$ is tame
and $|\bbX_0| = 4$, then $\bDelta$ is of type $(2, 2, 2, 2)$.

Fix $\lambda \in \bbX$. If $R_{\lambda, 0}$, \ldots, $R_{\lambda,
r_\lambda - 1}$ are chosen representatives of the isomorphisms
classes of the simple objects in $\add \calR_\lambda$, then we may
assume that $\tau R_{\lambda, i} = R_{\lambda, i - 1}$ for each $i
\in [0, r_\lambda - 1]$, where we put $R_{\lambda, i} :=
R_{\lambda, i \mod r_\lambda}$ for $i \in \bbZ$. For $i \in \bbZ$
and $n \in \bbN_+$ there exists a unique (up to isomorphism)
representation in $\calR_\lambda$, whose top and length in $\add
\calR_\lambda$ are $R_{\lambda, i}$ and $n$, respectively. We fix
such representation and denote it by $R_{\lambda, i}^{(n)}$, and
its dimension vector by $\be_{\lambda, i}^n$. Then the composition
factors of $R_{\lambda, i}^{(n)}$ are (starting from the top):
$R_{\lambda, i}$, $R_{\lambda, i - 1}$, \ldots, $R_{\lambda, i -
(n - 1)}$. Consequently, $\be_{\lambda, i}^n = \sum_{j \in [i - n
+ 1, i]} \be_{\lambda, j}$, where $\be_{\lambda, j} := \bdim
R_{\lambda, j}$ for $j \in \bbZ$. Moreover, if $i \in \bbZ$ and
$m, n \in \bbN_+$, then we have an exact sequence $0 \to
R_{\lambda, i - n}^{(m)} \to R_{\lambda, i}^{(m + n)} \to
R_{\lambda, i}^{(n)} \to 0$. Obviously, for each $R \in
\calR_\lambda$ there exist $i \in \bbZ$ and $n \in \bbN_+$ such
that $R \simeq R_{\lambda, i}^{(n)}$. Moreover, it is known that
the vectors $\be_{\lambda, 0}$, \ldots, $\be_{\lambda, r_\lambda -
1}$ are linearly independent. Consequently, if $R \in \add
\calR_\lambda$, then there exist uniquely determined $q_0^R,
\ldots, q_{r_\lambda - 1}^R \in \bbN$ such that $\bdim R = \sum_{i
\in [0, r_\lambda - 1]} q_i^R \be_{\lambda, i}$. Observe that the
numbers $q_{\lambda, 0}^R$, \ldots, $q_{\lambda, r_\lambda - 1}^R$
count the multiplicities in which the modules $R_{\lambda, 0}$,
\ldots, $R_{\lambda, r_\lambda - 1}$ appear as composition factors
in the Jordan--H\"older filtration of $R$ in the category $\add
\calR_\lambda$.

Let $R = \bigoplus_{\lambda \in \bbX} R_\lambda$ for $R_\lambda
\in \add \calR_\lambda$, $\lambda \in \bbX$. Then we put
$q_{\lambda, i}^R := q_i^{R_\lambda}$ for $\lambda \in \bbX$ and
$i \in [0, r_\lambda - 1]$. Next, we put $p_\lambda^R := \min \{
q_{\lambda, i}^R : \text{$i \in [0, r_\lambda - 1]$} \}$ for
$\lambda \in \bbX$, and $p_{\lambda, i}^R := q_{\lambda, i}^R -
p_\lambda^R$ for $\lambda \in \bbX$ and $i \in [0, r_\lambda -
1]$. Then
\[
\bdim R = \sum_{\lambda \in \bbX} p_\lambda^R \bh_\lambda +
\sum_{\lambda \in \bbX} \sum_{i \in [0, r_\lambda - 1]}
p_{\lambda, i}^R \be_{\lambda, i},
\]
where $\bh_\lambda := \sum_{i \in [0, r_\lambda - 1]}
\be_{\lambda, i}$ for $\lambda \in \bbX$. It is known that
$\bh_\lambda = \bh_\mu$ for any $\lambda, \mu \in \bbX$. We denote
this common value by $\bh$. Then
\[
\bdim R = p^R \bh + \sum_{\lambda \in \bbX} \sum_{i \in [0,
r_\lambda - 1]} p_{\lambda, i}^R \be_{\lambda, i},
\]
where $p^R := \sum_{\lambda \in \bbX} p_\lambda^R$. It is known
that $p^R = p^{R'}$ and $p_{\lambda, i}^R = p_{\lambda, i}^{R'}$
for any $\lambda \in \bbX$ and $i \in [0, r_\lambda - 1]$, if $R,
R' \in \add \calR$ and $\bdim R = \bdim R'$. Consequently, for
each $\bd \in \bR$ there exist uniquely determined $p^\bd \in
\bbN$ and $p_{\lambda, i}^{\bd} \in \bbN$, $\lambda \in \bbX$, $i
\in [0, r_\lambda - 1]$, such that for each $\lambda \in \bbX$
there exists $i \in [0, r_\lambda - 1]$ with $p_{\lambda, i}^{\bd}
= 0$ and
\[
\bd = p^{\bd} \bh + \sum_{\lambda \in \bbX} \sum_{i \in [0,
r_\lambda - 1]} p_{\lambda, i}^{\bd} \be_{\lambda, i}.
\]

Let $\lambda, \mu \in \bbX$, $i, j \in \bbZ$, and $m, n \in
\bbN_+$. Then
\[
\dim_k \Hom_{\bDelta} (R_{\lambda, i}^{(n)}, R_{\mu, j}^{(m)}) =
\min \{ q_{\lambda, i \mod r_\lambda}^{R_{\mu, j}^{(m)}}, q_{\mu,
(j - m + 1) \mod r_\lambda}^{R_{\lambda, i}^{(n)}} \}.
\]
In particular, if  $\lambda \in \bbX$, $i \in [0, r_\lambda - 1]$,
$n \in \bbN_+$, $R \in \add \calR$ and $\Hom_{\bDelta}
(R_{\lambda, i}^{(n)}, R) \neq 0$, then $q_{\lambda, i}^R \neq 0$.
Moreover, the above formula together with the Auslander--Reiten
formula imply that
\begin{gather*}
\langle \be_{i, \lambda}^n, \bd \rangle_{\bDelta} = p_{\lambda, i
\mod r_\lambda}^{\bd} - p_{\lambda, (i - n) \mod r_\lambda}^{\bd}
\\
\intertext{and} %
\langle \bd,  \be_{i, \lambda}^n \rangle_{\bDelta} = p_{\lambda,
(i - n + 1) \mod r_\lambda}^{\bd} - p_{\lambda, (i + 1) \mod
r_\lambda}^{\bd}
\end{gather*}
for any $\lambda \in \bbX$, $i \in \bbZ$, $n \in \bbN_+$, and $\bd
\in \bR$. Consequently, $\langle \bh, \bd \rangle_{\bDelta} = 0 =
\langle \bd, \bh \rangle_{\bDelta}$ for each $\bd \in \bR$. In
particular, $q_{\bDelta} (\bh) = 0$. On the other hand, if $\bd
\in \bR$, then $q_{\bDelta} (\bd) = 0$ if and only if $\bd =
p^{\bd} \bh$. One also shows that $\bh$ is indivisible.

We also need some other properties of the Tits form, which we list
now.

\begin{proposition} \label{proposition Tits}
Assume that $\bDelta$ is tame. Then the following hold.
\begin{enumerate}

\item \label{point Tits1}
$q_{\bDelta} (\bd) \geq 0$ for each dimension vector $\bd$.

\item \label{point Tits2}
If $q_{\bDelta} (\bd) = 0$ for a dimension vector $\bd$, then $\bd
\in \bP \cup \bR \cup \bQ$ and  $\langle \bd, \bd_0
\rangle_{\bDelta} + \langle \bd_0, \bd \rangle_{\bDelta} = 0$ for
each dimension vector $\bd_0$.

\item \label{point Tits3}
If $\bd \in \bP \cup \bQ$ is non-zero, then $\langle \bd, \bh
\rangle_{\bDelta} \neq 0$.

\item \label{point Tits4}
If $\bd \in \bP \cup \bQ$ is non-zero and $q_{\bDelta} (\bd) = 0$,
then $\langle \bd, \bd_0 \rangle_{\bDelta} \neq 0$ for each
non-zero vector $\bd_0 \in \bR$. In particular,
\[
|\langle \bd, \bh \rangle_{\bDelta}| \geq \max \{ r_\lambda :
\text{$\lambda \in \bbX$} \}.
\]

\item \label{point Tits5}
If there exists non-zero $\bd \in \bP \cup \bQ$ such that
$q_{\bDelta} (\bd) = 0$, then $\sum_{\lambda \in \bbX_0}
\frac{1}{r_\lambda} = |\bbX_0| - 2$. In particular, if this is the
case, then $\max \{ r_\lambda : \text{$\lambda \in \bbX$} \} \geq
2$ and $\max \{ r_\lambda : \text{$\lambda \in \bbX$} \} = 2$ if
and only if $\bDelta$ is of type $(2, 2, 2, 2)$.

\end{enumerate}
\end{proposition}

As a consequence we obtain the following.

\begin{corollary} \label{corollary inequality}
Let $\bd \in \bR$, $\bd' \in \bP + \bR$ and $\bd'' \in \bQ$. If
$p^{\bd} > 0$, $\bd' + \bd'' = \bd$ and $\bd'' \neq 0$, then
$\langle \bd'', \bd' \rangle_{\bDelta} \leq - p^{\bd} - 1$.
Moreover, $\langle \bd'', \bd' \rangle_{\bDelta} = - p^{\bd} - 1$
if and only if one of the following conditions is satisfied:
\begin{enumerate}

\item
$q_{\bDelta} (\bd'') = 1$ and $\langle \bd'', \bd
\rangle_{\bDelta} = -p^{\bd}$ \textup{(}in particular, $\langle
\bd'', \bh \rangle_{\bDelta} = -1$\textup{)}, or

\item
$q_{\bDelta} (\bd'') = 0$ and $\langle \bd'', \bd
\rangle_{\bDelta} = -2$.

\end{enumerate}
\end{corollary}

\begin{proof}
Put $\bd_0 := \bd - p^{\bd} \bh$. Then $\bd_0 \in \bR$. We have
\[
\langle \bd'', \bd' \rangle_{\bDelta} = \langle \bd'', \bd - \bd''
\rangle_{\bDelta} = - q_{\bDelta} (\bd'') + p^{\bd} \langle \bd'',
\bh \rangle_{\bDelta} + \langle \bd'', \bd_0 \rangle_{\bDelta}.
\]
Now $\langle \bd'', \bd_0 \rangle_{\bDelta} \leq 0$. Moreover,
$\langle \bd'', \bh \rangle_{\bDelta} \leq -1$ and $q_{\bDelta}
(\bd'') \geq 0$. Finally, if $q_{\bDelta} (\bd'') = 0$, then
$\langle \bd'', \bh \rangle_{\bDelta} \leq -2$, hence the
inequality follows.

These considerations also imply that $\langle \bd'', \bd'
\rangle_{\bDelta} = - p^{\bd} - 1$ if and only if one of the
following conditions is satisfied:
\begin{enumerate}

\item
$q_{\bDelta} (\bd'') = 1$, $\langle \bd'', \bh \rangle_{\bDelta} =
-1$ and $\langle \bd'', \bd_0 \rangle_{\bDelta} = 0$, or

\item
$q_{\bDelta} (\bd'') = 0$, $p^{\bd} = 1$, $\langle \bd'', \bh
\rangle_{\bDelta} = -2$ and $\langle \bd'', \bd_0
\rangle_{\bDelta} = 0$.

\end{enumerate}
These conditions immediately lead to (and follows from) the
conditions given in the corollary.
\end{proof}

We call a dimension vector $\bd \in \bR$ singular if $p^{\bd} > 0$
and there exists a dimension vector $\bx$ such that $\bx \leq
\bd$, $q_{\bDelta} (\bx) = 0$ and $|\langle \bx, \bd
\rangle_{\bDelta}| = 2$. It follows from the below proposition
that this definition coincides the the definition given in the
introduction.

\begin{proposition} \label{proposition singular}
Let $\bd \in \bR$ be such that $p^{\bd} > 0$.
\begin{enumerate}

\item \label{point singular1}
If $\bd$ is singular, then $\bd = \bh$ and $\bDelta$ is of type
$(2, 2, 2, 2)$.

\item \label{point singular2}
There exist $\bd' \in \bP + \bR$ and $\bd'' \in \bQ$ such that
$\bd' + \bd'' = \bd$, $q_{\bDelta} (\bd'') = 0$ and $\langle
\bd'', \bd \rangle_{\bDelta} = -2$, if and only if $\bd$ is
singular.

\end{enumerate}
\end{proposition}

\begin{proof}
\eqref{point singular1}~Fix a dimension vector $\bx$ such that
$\bx \leq \bd$, $q_{\bDelta} (\bx) = 0$ and $|\langle \bx, \bd
\rangle_{\bDelta}| = 2$. Proposition~\ref{proposition
Tits}\eqref{point Tits2} implies that $\bx \in \bP \cup \bR \cup
\bQ$. Since $\langle \bx, \bd \rangle_{\bDelta} \neq 0$, $\bx \not
\in \bR$. In particular, $\bx$ is non-zero. By symmetry, we may
assume $\bx \in \bP$. If $\bd_0 := \bd - p^{\bd} \bh$, then $2 =
p^{\bd} \langle \bx, \bh \rangle_{\bDelta} + \langle \bx, \bd_0
\rangle_{\bDelta}$. Using Proposition~\ref{proposition
Tits}\eqref{point Tits4} and~\eqref{point Tits5} we obtain that
$p^{\bd} = 1$ and $\bd_0 = 0$, i.e. $\bd = \bh$. Moreover,
$\bDelta$ must be of type $(2, 2, 2, 2)$ by
Proposition~\ref{proposition Tits}\eqref{point Tits5}.

\eqref{point singular2}~One implication is obvious. Now assume
there exists a dimension vector $\bx$ such that $\bx \leq \bd$,
$q_{\bDelta} (\bx) = 0$ and $|\langle \bx, \bd \rangle_{\bDelta}|
= 2$. From~\eqref{point singular1} we know that $\bd = \bh$. Easy
calculations show that $\langle \bh, \bh - \bx \rangle_{\bDelta} =
- \langle \bh, \bx \rangle_{\bDelta}$ and $q_{\bDelta} (\bh - \bx)
= 0$. Thus, Proposition~\ref{proposition Tits}\eqref{point Tits2}
implies that, up to symmetry, $\bx \in \bP$ and $\bh - \bx \in
\bQ$, and the claim follows.
\end{proof}

We finish this section with an example showing that singular
dimension vectors exist. Fix $\lambda \in k \setminus \{ 0, 1 \}$.
Let $\Delta$ be the quiver
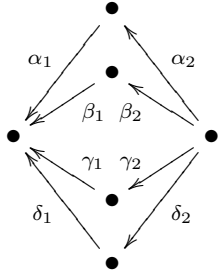
\[
\xymatrix@R=1eM{%
& \bullet \ar[ldd]_{\alpha_1}
\\
& \bullet \ar[ld]^(0.33){\beta_1}
\\
\bullet & & \bullet
\ar[luu]_{\alpha_2} \ar[lu]^(0.67){\beta_2} %
\ar[ld]_(0.67){\gamma_2} \ar[ldd]^{\delta_2}
\\
& \bullet \ar[lu]_(0.33){\gamma_1}
\\
& \bullet \ar[luu]^{\delta_1}
}%
\]
and $R := \{ \alpha_1 \alpha_2 + \beta_1 \beta_2 + \gamma_1
\gamma_2, \alpha_1 \alpha_2 + \beta_1 \beta_2 + \lambda \delta_1
\delta_2 \}$. Then $\bDelta$ is a concealed-canonical algebra of
type $(2, 2, 2, 2)$ (in fact, it is one of Ringel's canonical
algebras~\cite{Ringel1980}). Moreover, the vector $
\begin{smallmatrix}
& 2 \\ & 2 \\ 3 & & 1 \\ & 2 \\ & 2
\end{smallmatrix}
$ is singular -- the corresponding vector $\bx$ can be taken to be
$
\begin{smallmatrix}
& 1 \\ & 1 \\ 1 & & 1 \\ & 1 \\ & 1
\end{smallmatrix}
$ (the other choice is $
\begin{smallmatrix}
& 1 \\ & 1 \\ 2 & & 0 \\ & 1 \\ & 1
\end{smallmatrix}
$).

\section{Varieties of representations} \label{section varieties}

First we recall some facts from algebraic geometry. Let $\calX$ be
a closed subvariety of an affine space $\bbA^n$, $n \in \bbN$. We
say that $\calX$ is a complete intersection if there exist
polynomials $f_1, \ldots, f_m \in k [\bbA^n]$ such that $\dim
\calX = n - m$ and
\[
\{ f \in k [\bbA^n] : \text{$f (x) = 0$ for each $x \in \calX$} \}
= (f_1, \ldots, f_m).
\]
For $x \in \calX$ we denote by $T_x \calX$ the tangent space to
$\calX$ at $x$. We will use the following consequences of Serre's
criterion (see for example~\cite{Eisenbud1995}*{Theorem~18.15}).

\begin{proposition} \label{proposition Serre}
Let $\calX$ be a complete intersection.
\begin{enumerate}

\item \label{point normal}
Let $\calU := \{ x \in \calX : \text{$\dim_k T_x \calX = \dim
\calX$} \}$. Then $\calX$ is normal if and only if $\dim (\calX
\setminus \calU) < \dim \calX - 1$.

\item \label{point reduced}
Let $f_1, \ldots, f_m \in k [\calX]$,
\begin{gather*}
\calY := \{ x \in \calX : \text{$f_i (x) = 0$ for each $i \in [1,
m]$} \}
\\
\intertext{and} %
\calU := \{ x \in \calY : \text{$\partial f_1 (x)$, \ldots,
$\partial f_m (x)$ are linearly independent} \}.
\end{gather*}
If $\calU \cap \calC \neq \varnothing$ for each irreducible
component $\calC$ of $\calY$, then
\[
\{ f \in k [\calX] : \text{$f (x) = 0$ for each $x \in \calY$} \}
= (f_1, \ldots, f_m).
\]
In particular, $\calY$ is a complete intersection of dimension
$\dim \calX - m$. \qed

\end{enumerate}
\end{proposition}

Let $\Delta$ be a bound quiver and $\bd$ a dimension vector. By
$\rep_\Delta (\bd)$ we denote the set of the representations $M$
of $\Delta$ such that $M (x) = k^{\bd (x)}$ for each $x \in
\Delta_0$. We may identify $\rep_\bDelta (\bd)$ with the affine
space $\prod_{\alpha \in \Delta_1} \bbM_{\bd (t \alpha) \times \bd
(s \alpha)} (k)$. The group $\GL (\bd) := \prod_{x \in \Delta_0}
\GL (\bd (x))$ acts on $\rep_\Delta (\bd)$ by conjugation: $(g
\cdot M) (\alpha) := g (t \alpha) M (\alpha) g (s \alpha)^{-1}$
for $g \in \GL (\bd)$, $M \in \rep_\Delta(\bd)$ and $\alpha \in
\Delta_1$. Under this action the $\GL (\bd)$-orbits in
$\rep_\Delta (\bd)$ correspond to the isomorphism classes of the
representations of $\Delta$ with dimension vector $\bd$. We denote
the $\GL (\bd)$-orbit of a representation $M \in \rep_\Delta
(\bd)$ by $\calO (M)$.

Now let $\bDelta$ be a bound quiver and $\bd$ a dimension vector.
By $\rep_\bDelta (\bd)$ we denote the intersection of $\rep_\Delta
(\bd)$ with $\rep \bDelta$. Then $\rep_{\bDelta} (\bd)$ is a
closed $\GL (\bd)$-invariant subset of $\rep_\Delta (\bd)$ and we
call it the variety of representations of $\bDelta$ of dimension
vector $\bd$. If $M, N \in \rep_{\bDelta} (\bd)$ and there exists
an exact sequence $0 \to N' \to M \to N'' \to 0$ such that $N
\simeq N' \oplus N''$, then $N \in \ol{\calO (M)}$. In particular,
$S^{\bd} \in \ol{\calO (M)}$ for each $M \in \rep_\Delta (\bd)$.
If $\calV$ is a $\GL (\bd)$-invariant subset of $\rep_\Delta
(\bd)$ and $M \in \calV$, then we say that the orbit $\calO (M)$
is maximal in $\calV$ if $\calO (N) = \calO (M)$ for each $N \in
\calV$ such that $\calO (M) \subseteq \ol{\calO (N)}$.

Put $a_{\bDelta} (\bd) := \dim \GL (\bd) - q_{\bDelta} (\bd)$ for
a bound quiver $\bDelta$ and a dimension vector $\bd$. The
following facts were proved in~\cite{BobinskiSkowronski2002}.

\begin{proposition} \label{proposition variety}
Let $\bd$ be the dimension vector of a periodic representation
over a tame concealed-canonical bound quiver $\bDelta$. Then the
following hold.
\begin{enumerate}

\item \label{point variety 1}
The variety $\rep_{\bDelta} (\bd)$ is a normal complete
intersection of dimension $a_{\bDelta} (\bd)$.

\item \label{point variety 3}
If there exists $M \in \rep_{\bDelta} (\bd)$ such that
$\Ext_{\bDelta}^1 (M, M) = 0$, then $\ol{\calO (M)} =
\rep_{\bDelta} (\bd)$.

\item \label{point variety 2}
If $\Ext_{\bDelta}^1 (M, M) \neq 0$ for each $M \in \rep_{\bDelta}
(\bd)$, then there exists a convex subquiver $\bDelta'$ of
$\bDelta$ and a sincere separating exact subcategory $\calR'$ in
$\rep \bDelta'$ such that $M \in \add \calR'$ for each maximal
orbit $\calO (M)$ in $\rep_{\bDelta} (\bd)$.

\item \label{point variety 4}
If $M \in \rep_{\bDelta} (\bd)$, then there is a canonical
epimorphism
\[
\pi_M : T_M \rep_{\bDelta} (\bd) \to \Ext_{\bDelta}^1 (M, M)
\]
with kernel $T_M \calO (M)$. \qed

\end{enumerate}
\end{proposition}

Let $\bd$ be the dimension vector of a periodic module over a tame
concealed-canonical bound quiver $\bDelta$. The above theorem
implies that in order to prove that $\ol{\calO (M)}$ is a normal
complete intersection for each maximal orbit $\calO (M)$ in
$\rep_{\bDelta} (\bd)$, we may assume that $\bd$ is the dimension
vector of a direct sum of modules from a sincere separating exact
subcategory of $\ind \bDelta$. Thus we fix a tame bound quiver
$\bDelta$ and a sincere separating exact subcategory $\calR$ of
$\ind \bDelta$. We will use freely notation introduced in
Section~\ref{section tubular}. It follows
from~\cite{BobinskiSkowronski2002}*{Section~3} that if $\bd \in
\bR$, then $M \in \add \calR$ for each maximal orbit $\calO (M)$
in $\rep_{\bDelta} (\bd)$.

For a full subcategory $\calX$ of $\ind \bDelta$ and a dimension
vector $\bd$ we denote by $\calX (\bd)$ the intersection of
$\rep_{\bDelta} (\bd)$ with $\add \calX$. If $\bd', \bd' \in
\bbN^{\Delta_0}$, $C' \subseteq \rep_{\bDelta} (\bd')$ and $C''
\subseteq \rep_{\bDelta} (\bd'')$, then we denote by $C' \oplus
C''$ the subset of $\rep_{\bDelta} (\bd' + \bd'')$ consisting of
all $M$ such that $M \simeq M' \oplus M''$ for some $M' \in C'$
and $M'' \in C''$. The following fact follows
from~\cite{Bobinski2007}*{Section~3}.

\begin{proposition} \label{proposition oplus}
If $\bd' \in \bP + \bR$ and $\bd'' \in \bQ$, then $(\calP \cup
\calR) (\bd') \oplus \calQ (\bd'')$ is an irreducible
constructible subset of $\rep_{\bDelta} (\bd' + \bd'')$ of
dimension $a_{\bDelta} (\bd) + \langle \bd'', \bd'
\rangle_{\bDelta}$. \qed
\end{proposition}

Using Corollary~\ref{corollary inequality} we immediately get the
following.

\begin{corollary} \label{corollary dimension}
Let $\bd \in \bR$, $\bd' \in \bP + \bR$ and $\bd'' \in \bQ$. If
$p^{\bd} > 0$, $\bd' + \bd'' = \bd$ and $\bd'' \neq 0$, then
\[
\dim ((\calP \cup \calR) (\bd') \oplus \calQ (\bd'')) \leq
a_{\bDelta} (\bd) - p^{\bd} - 1.
\]
Moreover, the equality holds if and only if one of the following
conditions is satisfied:
\begin{enumerate}

\item
$q_{\bDelta} (\bd'') = 1$ and $\langle \bd'', \bd
\rangle_{\bDelta} = -p^{\bd}$ \textup{(}in particular, $\langle
\bd'', \bh \rangle_{\bDelta} = -1$\textup{)}, or

\item
$q_{\bDelta} (\bd'') = 0$ and $\langle \bd'', \bd
\rangle_{\bDelta} = -2$ \textup{(}in particular, $\bDelta$ is of
type $(2, 2, 2, 2)$ and  $\bd = \bh$\textup{)}. \qed

\end{enumerate}
\end{corollary}

Observe that
\[
\rep_{\bDelta} (\bd) = \calR (\bd) \cup \bigcup_{\substack{\bd'
\in \bP + \bR, \; \bd'' \in \bQ \\ \bd' + \bd'' = \bd, \; \bd''
\neq 0}} (\calP \cup \calR) (\bd') \oplus \calQ (\bd'')
\]
for each $\bd \in \bR$. Indeed, if $M \in (\calP \cup \calR)
(\bd)$ and we write $M = M' \oplus M''$ for $M' \in \add \calP$
and $M'' \in \add \calR$, then $\langle \bdim M', \bh
\rangle_{\bDelta} = \langle \bd, \bh \rangle_{\bDelta} = 0$, hence
$M' = 0$ by Proposition~\ref{proposition Tits}\eqref{point Tits3}.
The above formula together with Corollary~\ref{corollary
dimension} implies that $\dim (\rep_{\bDelta} (\bd) \setminus
\calR (\bd)) \leq a_{\bDelta} (\bd) - p^{\bd} - 1$.

\section{Stability and semi-invariants} \label{section
semiinvariants}

Let $\Delta$ be a quiver and $\theta \in \bbZ^{\Delta_0}$. We
treat $\theta$ as a $\bbZ$-linear function $\bbZ^{\Delta_0} \to
\bbZ$ in a usual way. A representation $M$ of $\Delta$ is called
$\theta$-semi-stable if $\theta (\bdim M) = 0$ and $\theta (\bdim
N) \geq 0$ for each subrepresentation $N$ of $M$. The full
subcategory of $\theta$-semi-stable representations of $\Delta$ is
an exact subcategory of $\rep \Delta$. Two $\theta$-semi-stable
representations are called S-equivalent if they have the same
composition factors within this category. If $\bd$ is a dimension
vector, then by a semi-invariant of weight $\theta$ we mean every
function $c \in k [\rep_\Delta (\bd)]$ such that $c (g \cdot M) =
\chi^\theta (g) c (M)$ for any $g \in \GL (\bd)$ and $M \in
\rep_\Delta (\bd)$, where $\chi^\theta (g) := \prod_{x \in
\Delta_0} (\det g (x))^{\theta (x)}$ for $g \in \GL (\bd)$.

Now let $\bDelta$ be a bound quiver and $\bd$ a dimension vector.
If $\theta \in \bbZ^{\Delta_0}$, then a function $c \in k
[\rep_{\bDelta} (\bd)]$ is called a semi-invariant of weight
$\theta$ if $c$ is a restriction of a semi-invariant of weight
$\theta$ from $k [\rep_\Delta (\bd)]$. This definition differs
from the definition used in other papers on the subject (see for
example~\cites{BobinskiRiedtmannSkowronski2008, DerksenWeyman2002,
Domokos2002, DomokosLenzing2002}), however it is consistent with
King's approach~\cite{King1994}. We denote the space of the
semi-invariants of weight $\theta$ by $\SI [\bDelta, \bd]_\theta$.
If $\theta \in \bbZ^{\Delta_0}$, then we put $\Lambda_\theta
(\bd):= \bigoplus_{n \in \bbN} \SI [\bDelta, \bd]_{n \theta}$.
Note that $\Lambda_\theta (\bd)$ is a graded ring. For $M \in
\rep_{\bDelta} (\bd)$ we denote by $\calI_\theta (M)$ the ideal in
$\Lambda_\theta (\bd)$ generated by the homogeneous elements $c$
such that $c (M) = 0$.

The following results were proved in~\cite{King1994}.

\begin{proposition} \label{proposition King}
Let $\bDelta$ be a bound quiver, $\bd$ a dimension vector, and
$\theta \in \bbZ^{\Delta_0}$.
\begin{enumerate}

\item
If $M \in \rep_{\bDelta} (\bd)$, then $M$ is $\theta$-semi-stable
if and only if there exists a semi-invariant $c$ of weight $n
\theta$ for some $n \in \bbN_+$ such that $c (M) \neq 0$.

\item \label{point King2}
If $M, N \in \rep_{\bDelta} (\bd)$ are $\theta$-semi-stable, then
$M$ and $N$ are S-equivalent if and only if $\calI_\theta (M) =
\calI_\theta (N)$. \qed

\end{enumerate}
\end{proposition}

Now we recall a construction from~\cite{Domokos2002}. Let
$\bDelta$ be a bound quiver. Fix a representation $V$ of
$\bDelta$. We define $\theta^V : \bbZ^{\Delta_0} \to \bbZ$ by the
condition
\[
\theta^V (\bdim M) = - \dim_k \Hom_{\bDelta} (V, M) + \dim_k
\Hom_{\bDelta} (M, \tau V)
\]
for each representation $M$ of $\bDelta$. The Auslander--Reiten
formula implies that $\theta^V = - \langle \bdim V, - \rangle$ if
$\pdim_{\bDelta} V \leq 1$. Dually, if $\idim_{\bDelta} V \leq 1$,
then $\theta^V = \langle -, \bdim \tau V \rangle$.

Now let $\bd$ be a dimension vector. If $\theta^V (\bd) = 0$, then
we define a function $c^V \in k [\rep_{\bDelta} (\bd)]$ in the
following way. Let $P_1 \xrightarrow{f} P_0 \to V \to 0$ be a
minimal projective presentation of $V$. There exist vertices
$x_1$, \ldots, $x_n$, $y_1$, \ldots, $y_m$ of $\Delta$ such that
$P_1 = \bigoplus_{i \in [1, n]} P_{x_i}$ and $P_0 = \bigoplus_{j
\in [1, m]} P_{y_j}$. Moreover, there exist $\omega_{i, j} \in k
\bDelta (y_j, x_i)$, $i \in [1, n]$, $j \in [1, m]$, such that $f
=
\begin{bmatrix}
P_{\omega_{i, j}}
\end{bmatrix}_{\substack{j \in [1, m] \\ i \in [1, n]}}$.
Consequently, if $M \in \rep_{\bDelta}( \bd)$, then
\[
\Hom_{\bDelta} (f, M) =
\begin{bmatrix}
M (\omega_{i, j})
\end{bmatrix}_{\substack{i \in [1, n] \\ j \in [1, m]}} :
\bigoplus_{j \in [1, m]} M (y_j) \to \bigoplus_{i \in [1, n]} M
(x_i).
\]
In addition, one shows $\dim_k \Ker \Hom_{\bDelta} (f, M) = \dim_k
\Hom_{\bDelta} (V, M)$ and $\dim_k \Coker \Hom_{\bDelta} (f, M) =
\dim_k \Hom_{\bDelta} (M, \tau V)$. Consequently,
\begin{multline*}
\sum_{j \in [1, m]} \dim_k M (y_j) - \sum_{i \in [1, n]} \dim_k M
(x_i)
\\
= \dim_k \Hom_{\bDelta} (M, \tau V) - \dim_k \Hom_{\bDelta} (V, M)
= \theta^V (\bd) = 0.
\end{multline*}
Thus, it makes sense to define $c^V \in k [\rep_{\bDelta} (\bd)]$
by
\[
c^V (M) := \det \Hom_{\bDelta} (f, M)
\]
for $M \in \rep_{\bDelta} (\bd)$. Note that $c^V (M) = 0$ if and
only if $\Hom_{\bDelta} (V, M) \neq 0$. It is known that $c^V \in
\SI [\bDelta, \bd]_{\theta^V}$. This function depends on the
choice of $f$, but functions obtained for different $f$'s differ
only by non-zero scalars. In fact, we could start with an
arbitrary projective presentation $P_1 \xrightarrow{f} P_0 \to V
\to 0$ of $V$ such that $\dim_k \Hom_{\bDelta} (P_1, M) = \dim_k
\Hom_{\bDelta} (P_0, M)$. As an easy consequence we obtain the
following (see~\cite{DerksenWeyman2002}*{Proposition~2}
and~\cite{Domokos2002}*{Lemma~3.3}).

\begin{lemma} \label{lemma multsemi}
Let $\bDelta$ be a bound quiver and $\bd$ a dimension vector.
\begin{enumerate}

\item \label{point multsemi1}
If $V = V_1 \oplus V_2$, $\theta^V (\bd) = 0$ and $c^V \neq
0$, then $\theta^{V_1} (\bd) = 0 = \theta^{V_2} (\bd)$.

\item \label{point multsemi2}
If $0 \to V_1 \to V \to V_2 \to 0$ and $\theta^V (\bd) =
\theta^{V_1} (\bd) = \theta^{V_2} (\bd) = 0$, then $c^V = c^{V_1}
c^{V_2}$. \qed

\end{enumerate}
\end{lemma}

The following result follows from the proof
of~\cite{Domokos2002}*{Theorem~3.2} (note that the assumption
about the characteristic of $k$ made
in~\cite{Domokos2002}*{Theorem~3.2} is only necessary for
surjectivity of the restriction morphism, which we have for free
with our definition of semi-invariants).

\begin{proposition}
Let $\bDelta$ be a bound quiver and $\bd$ a dimension vector. If
$\theta \in \bbZ^{\Delta_0}$, then the space $\SI [\bDelta,
\bd]_{\theta}$ is spanned by the functions $c^V$ for $V \in \rep
\bDelta$ such that $\theta^V = \theta$. \qed
\end{proposition}

Now we apply our considerations in the case of tame
concealed-canonical quivers. For the rest of the section we fix a
tame bound quiver $\bDelta$ and a sincere separating exact
subcategory $\calR$ of $\ind \bDelta$. We will use notation
introduced in Section~\ref{section tubular}. We fix $\bd \in \bR$
such that $p^{\bd} > 0$ and put $\theta := - \langle \bh, -
\rangle_{\bDelta}$.

First observe that $M \in \rep \bDelta$ is $\theta$-semi-stable if
and only if $M \in \add \calR$. Consequently, if $M$ and $N$ are
$\theta$-semi-stable, then $M$ and $N$ are S-equivalent if and
only if $q_{\lambda, i}^M = q_{\lambda, i}^N$ for any $\lambda \in
\bbX$ and $i \in [0, r_\lambda - 1]$. In particular, there are
only finitely many isomorphism classes in each S-equivalence
class.

Now fix $V \in \rep \bDelta$ such that $\theta^V = n \theta$ for
some $n \in \bbN$ and $c^V \neq 0$. We show that $V \in \add
\calR$ and $\bdim V = n \bh$. Indeed, write $V = P \oplus R \oplus
Q$ for $P \in \add \calP$, $R \in \add \calR$ and $Q \in \add
\calQ$. If $P \neq 0$, then $\theta^P (\bd) \leq - \langle \bdim
P, \bh \rangle_{\bDelta} < 0$ by Proposition~\ref{proposition
Tits}\eqref{point Tits3}, hence $c^V = 0$ by Lemma~\ref{lemma
multsemi}\eqref{point multsemi1}. Consequently, $P = 0$ and,
dually, $Q = 0$, thus $V = R \in \add \calR$. In particular,
$\pdim_{\bDelta} V = 1$, hence $- \langle n \bh, - \rangle =
\theta^V = - \langle \bdim V, - \rangle_{\bDelta}$, and this
implies that $\bdim V = n \bh$.

For $\lambda \in \bbX$ we denote by $\calA_\lambda (\bd)$ the set
of all $i \in [0, r_\lambda - 1]$ such that $p_{\lambda, i}^{\bd}
= 0$. Next, for $i \in \calA_\lambda (\bd)$ we denote by
$n_{\lambda, i}$ the minimal $n \in \bbN_+$ such that $p_{\lambda,
(i - n) \mod r_\lambda}^{\bd} = 0$, and put $V_{\lambda, i} :=
R_{\lambda, i}^{(n_{\lambda, i})}$. Observe that
$\theta^{V_{\lambda, i}} (\bd) = - \langle \bdim V_{\lambda, i},
\bd \rangle_{\bDelta} = 0$ for any $\lambda \in \bbX$ and $i \in
\calA_\lambda (\bd)$. We put $c_{\lambda, i} := c^{V_{\lambda,
i}}$ for $\lambda \in \bbX$ and $i \in \calA_\lambda (\bd)$. More
generally, if $\lambda \in \bbX$ and $J \subseteq \calA_\lambda
(\bd)$, then we put $V_{\lambda, J} := \bigoplus_{i \in J}
V_{\lambda, i}$ and $c_{\lambda, J} := c^{V_{\lambda, J}} =
\prod_{i \in J} c_{\lambda, i}$. In particular, we put $V_\lambda
:= V_{\lambda, \calA_\lambda (\bd)}$ and $c_\lambda := c_{\lambda,
\calA_\lambda (\bd)}$ for $\lambda \in \bbX$. Then $c_\lambda \in
\SI [\bDelta, \bd]_\theta$ for each $\lambda \in \bbX$. Observe
that Lemma~\ref{lemma multsemi}\eqref{point multsemi2} implies
that $c_\lambda = c^{R_{\lambda, i}^{(r_\lambda)}}$ for any
$\lambda \in \bbX$ and $i \in \calA_\lambda (\bd)$. More general,
$c_\lambda^p = c^{R_{\lambda, i}^{(p r_\lambda)}}$ for any $p \in
\bbN_+$, $\lambda \in \bbX$ and $i \in \calA_\lambda (\bd)$.

We have the following information about $\Lambda_\theta (\bd)$.

\begin{proposition} \label{proposition Lambda}
Let $\bd \in \bR$ be such that $p^{\bd} > 0$ and $\theta := -
\langle \bh, - \rangle$. Then $\Lambda_\theta (\bd)$ is generated
by the functions $c_\lambda$, $\lambda \in \bbX$.
\end{proposition}

\begin{proof}
First we show that if $\lambda \in \bbX$, $i \in \bbZ$, $n \in
\bbN_+$, $\theta^{R_{\lambda, i}^{(n)}} (\bd) = 0$ and
$c^{R_{\lambda, i}^{(n)}} \neq 0$, then $p_{\lambda, i \mod
r_\lambda}^{\bd} = p_{\lambda, (i - n) \mod r_\lambda}^{\bd}$ and
$p_{\lambda, j \mod r_\lambda}^{\bd} \geq p_{\lambda, i \mod
r_\lambda}^{\bd}$ for each $j \in [i - n + 1, i - 1]$. Indeed, the
former condition follows from the equality $\langle \be_{\lambda,
i}^{(n)}, \bd \rangle =  - \theta^{R_{\lambda, i}^{(n)}} (\bd) =
0$. Moreover, if there exists $j \in [i - n + 1, i - 1]$ such that
$p_{\lambda, j \mod r_\lambda}^{\bd} < p_{\lambda, i \mod
r_\lambda}^{\bd}$, then $\Hom_{\bDelta} (R_{\lambda, i}^{(n)}, R)
\neq 0$ for each $R \in \calR (\bd)$, hence $c^{R_{\lambda,
i}^{(n)}} = 0$.

We have the following important consequence of the above
observation. Assume that $\lambda \in \bbX$, $i \in [0, r_\lambda
- 1]$, $p \in \bbN_+$ and $c^{R_{\lambda, i}^{(p r_\lambda)}} \neq
0$. Then $p_{\lambda, j }^{\bd} \geq p_{\lambda, i}^{\bd}$ for
each $j \in [0, r_\lambda - 1]$, hence $i \in \calA_\lambda
(\bd)$. In particular, $c^{R_{\lambda, i}^{(p r_\lambda)}} =
c_\lambda^p$.

Now assume that $R \in \rep \bDelta$, $\theta^R = n \theta$ for
some $n \in \bbN$, and $c^R \neq 0$. We know that $R \in \add
\calR$ and $\bdim R = n \bh$. If $R = \bigoplus_{\lambda \in \bbX}
R_\lambda$ for $R_\lambda \in \calR_\lambda$, $\lambda \in \bbX$,
then $\bdim R_\lambda = p_\lambda^R \bh$ for each $\lambda \in
\bbX$. We show that $c^{R_\lambda} = c_\lambda^{p_\lambda^R}$ for
each $\lambda \in \bbX$, hence the claim will follow from
Lemma~\ref{lemma multsemi}\eqref{point multsemi1}.

Fix $\lambda \in \bbX$ and write $R_\lambda = \bigoplus_{j \in [1,
m]} R_{\lambda, i_j}^{(n_j)}$ for $m \in \bbN_+$, $i_1, \ldots,
i_m \in \bbZ$ and $n_1, \ldots, n_m \in \bbN_+$. If $n_j \equiv 0
\pmod{r_\lambda}$ for each $j \in [1, m]$, then the claim follows.
Thus assume $n_1 \not \equiv 0 \pmod{r_\lambda}$. Since $\bdim
R_\lambda = p_\lambda^R \bh$, we may assume that $i_2 = i_1 -
n_1$. Then we have an exact sequence $0 \to R_{\lambda,
i_2}^{(n_2)} \to R_{\lambda, i_1}^{(n_1 + n_2)} \to R_{\lambda,
i_1}^{(n_1)} \to 0$, hence Lemma~\ref{lemma multsemi}\eqref{point
multsemi2} implies that $c^R = c^{R'}$, where $R' := R_{\lambda,
i_1}^{(n_1 + n_2)} \oplus \bigoplus_{j \in [3, n]} R_{\lambda,
i_j}^{(n_j)}$. Now the claim follows by induction.
\end{proof}

As a consequence we get the following.

\begin{corollary} \label{corollary Sequivalence}
Let $\bd \in \bR$ be such that $p^{\bd} > 0$ and $\theta := -
\langle \bh, - \rangle$. If $M, N \in \calR (\bd)$, then $M$ and
$N$ are S-equivalent if and only if there exists $\mu \in k$ such
that $c_\lambda (M) = \mu c_\lambda (N)$ for each $\lambda \in
\bbX$.
\end{corollary}

\begin{proof}
Follows immediately form Propositions~\ref{proposition
King}\eqref{point King2} and~\ref{proposition Lambda}.
\end{proof}

We list some consequences of the description of the maximal orbits
in $\rep_{\bDelta} (\bd)$ given
in~\cite{BobinskiSkowronski2002}*{Proposition~5} (see
also~\cite{Ringel1980}*{Theorem~3.5}). Recall that $M \in \calR
(\bd)$ for each maximal orbit $\calO (M)$ in $\rep_{\bDelta}
(\bd)$. Next, if $\calO (M)$ is maximal in $\rep_{\bDelta} (\bd)$,
then $\dim \calO (M) = a_{\bDelta} (\bd) - p^{\bd}$. In
particular, the maximal orbits in $\rep_{\bDelta} (\bd)$ coincide
with the orbits of maximal dimension. Moreover, if $\lambda \in
\bbX$, then there exists at most one $i \in \calA_\lambda (\bd)$
such that $c_{\lambda, i} (M) = 0$. We put
\[
\hat{\bbX} (M) := \{ (\lambda, i) : \text{$\lambda \in \bbX$, $i
\in \calA_\lambda (\bd)$ and $c_{\lambda, i} (M) = 0$} \},
\]
and denote by $\bbX (M)$ the image of $\hat{\bbX} (M)$ under the
projection on the first coordinate. If $\lambda \in \bbX$, then
$\lambda \in \bbX (M)$ if and only if $p_\lambda^M \neq 0$. In
particular, $|\bbX (M)| \leq p^{\bd}$. Finally, if $M, N \in
\rep_{\bDelta} (\bd)$ are S-equivalent, the orbits $\calO (M)$ and
$\calO (N)$ are maximal, and $\hat{\bbX} (M) \subseteq \hat{\bbX}
(N)$, then $\calO (M) = \calO (N)$.

For a representation $V$ of $\bDelta$ such that $\theta^V (\bd) =
0$ we denote by $\calH^V (\bd)$ the zero set of $c^V$, i.e.\
$\calH^V (\bd) := \{ M \in \rep_{\bDelta} (\bd) :
\text{$\Hom_{\bDelta} (V, M) \neq 0$} \}$. Moreover, we say that
an exact sequence $0 \to M \to N \to L \to 0$ is $V$-exact if the
induced sequence
\[
0 \to \Hom_{\bDelta} (V, M) \to \Hom_{\bDelta} (V, M) \to
\Hom_{\bDelta} (V, L) \to 0
\]
is exact. We need the following version
of~\cite{RiedtmannZwara2011}*{Corollary~7.4}.

\begin{proposition} \label{proposition Zwara}
Let $V$ be a representation of $\bDelta$ such that $\theta^V (\bd)
= 0$.
\begin{enumerate}

\item \label{point Zwara1}
If $M \in \calH^V (\bd)$ and $\dim_k \Hom_{\bDelta} (V, M) = 1$,
then
\[
\Ker \partial c^V (M) = \{ Z \in T_M \rep_{\bDelta} (\bd) :
\text{$\pi_M (Z)$ is $V$-exact} \}.
\]

\item \label{point Zwara2}
If $M \in \calH^V (\bd)$ and $\dim_k \Hom_{\bDelta} (V, M) \geq
2$, then $\Ker \partial c^V (M) = T_M \rep_{\bDelta} (\bd)$.

\end{enumerate}
\end{proposition}

\section{Auxiliary lemmas} \label{section lemmas}

Throughout this section we fix a tame bound quiver $\bDelta$ and a
sincere separating exact subcategory $\calR$ of $\ind \bDelta$. We
use freely notation introduced in Section~\ref{section tubular}.
We also fix $\bd \in \bR$ such that $p := p^{\bd} > 0$.

\begin{lemma} \label{lemma intersection}
If $\lambda_0, \ldots, \lambda_p \in \bbX$ are pairwise different,
then
\[
\bigcap_{l \in [0, p]} \calH^{V_{\lambda_l}} (\bd) =
\bigcap_{\lambda \in \bbX}  \calH^{V_\lambda} (\bd) =
\bigcup_{\substack{\bd' \in \bP + \bR, \; \bd'' \in \bQ \\ \bd' +
\bd'' = \bd, \; \bd'' \neq 0}} (\calP \cup \calR) (\bd') \oplus
\calQ (\bd'').
\]
\end{lemma}

\begin{proof}
Obviously, $\bigcap_{l \in [0, p]} \calH^{V_{\lambda_l}} (\bd)
\supseteq \bigcap_{\lambda \in \bbX}  \calH^{V_\lambda} (\bd)$.

Now fix $\lambda \in \bbX$, $\bd' \in \bP + \bR$ and $\bd'' \in
\bQ$ such that $\bd' + \bd'' = \bd$ and $\bd'' \neq 0$. If $P \in
\calP (\bd')$ and $Q \in \calQ (\bd'')$, then
Proposition~\ref{proposition Tits}\eqref{point Tits3} implies that
\[
\dim_k \Hom_{\bDelta} (V_\lambda, P \oplus Q) \geq \dim_k
\Hom_{\bDelta} (V_\lambda, Q) = \langle \bh, \bd''
\rangle_{\bDelta} > 0,
\]
hence $(\calP \cup \calR) (\bd') \oplus \calQ (\bd'') \subseteq
\calH^{V_\lambda} (\bd)$.

Finally, assume that $R \in \calR (\bd) \cap \bigcap_{l \in [0,
p]} \calH^{V_{\lambda_l}} (\bd)$. Then $p_{\lambda_l}^R > 0$ for
each $l \in [0, p]$. Consequently, $p^R \geq \sum_{l \in [0, p]}
p_{\lambda_l}^R > p$, a contradiction.
\end{proof}

\begin{corollary} \label{corollary intersection}
Let $\lambda_0, \ldots, \lambda_p \in \bbX$ be pairwise different.
If $\calC$ is an irreducible component of $\bigcap_{l \in [0, p]}
\calH^{V_{\lambda_l}} (\bd)$, then $\dim \calC = a_{\bDelta} (\bd)
- p - 1$ and there exist $\bd' \in \bP + \bR$ and $\bd'' \in \bQ$
such that $\bd' + \bd'' = \bd$, $\bd'' \neq 0$ and $\calC =
\ol{(\calP \cup \calR) (\bd') \oplus \calQ (\bd'')}$.
\end{corollary}

\begin{proof}
It follows from Lemma~\ref{lemma intersection} that $\calC$ is an
irreducible component of $\ol{(\calP \cup \calR) (\bd') \oplus
\calQ (\bd'')}$ for some $\bd' \in \bP + \bR$ and $\bd'' \in \bQ$
such that $\bd' + \bd'' = \bd$ and $\bd'' \neq 0$. Since $(\calP
\cup \calR) (\bd') \oplus \calQ (\bd'')$ is irreducible by
Proposition~\ref{proposition oplus}, $\calC = \ol{(\calP \cup
\calR) (\bd') \oplus \calQ (\bd'')}$.

We know from Proposition~\ref{proposition variety}\eqref{point
variety 1} that $\dim \rep_{\bDelta} (\bd) = a_{\bDelta} (\bd)$,
hence Krull's Principal Ideal
Theorem~\cite{Kunz1985}*{Section~V.3} implies that $\dim \calC
\geq a_{\bDelta} (\bd) - p - 1$. On the other hand, $\dim \calC =
\dim (\calP \cup \calR) (\bd') \oplus \calQ (\bd'') \leq
a_{\bDelta} (\bd) - p - 1$ by Corollary~\ref{corollary dimension},
and the claim follows.
\end{proof}

\begin{lemma} \label{lemma intersectionbis}
Let $\lambda_0, \ldots, \lambda_p \in \bbX$ be pairwise different
and $J_l \subseteq \calA_{\lambda_l} (\bd)$, $l \in [0, p]$. If
$\calC$ is an irreducible component of $\bigcap_{l \in [0, p]}
\calH^{V_{\lambda_l, J_l}} (\bd)$, then $\calC$ is an irreducible
component of $\bigcap_{l \in [0, p]} \calH^{V_{\lambda_l}} (\bd)$.
\end{lemma}

\begin{proof}
Similarly as in the proof of Corollary~\ref{corollary
intersection} we show that $\dim \calC \geq a_{\bDelta} (\bd) - p
- 1$. On the other hand, $\calC \subseteq \bigcap_{l \in [0, p]}
\calH^{V_{\lambda_l}} (\bd)$, hence there exists an irreducible
component $\calC'$ of $\bigcap_{l \in [0, p]}
\calH^{V_{\lambda_l}} (\bd)$ such that $\calC \subseteq \calC'$.
Corollary~\ref{corollary intersection} says that $\dim \calC' =
a_{\bDelta} (\bd) - p - 1$, hence $\calC = \calC'$.
\end{proof}

\begin{corollary} \label{corollary intersectionbis}
Let $\lambda_0, \ldots, \lambda_p \in \bbX$ be pairwise different
and $J_l \subseteq \calA_{\lambda_l} (\bd)$, $l \in [0, p]$. If
$\calC$ is an irreducible component of $\bigcap_{l \in [0, p]}
\calH^{V_{\lambda_l, J_l}} (\bd)$, then $\dim \calC = a_{\bDelta}
(\bd) - p - 1$ and there exist $\bd' \in \bP + \bR$ and $\bd'' \in
\bQ$ such that $\bd' + \bd'' = \bd$, $\bd'' \neq 0$ and $\calC =
\ol{(\calP \cup \calR) (\bd') \oplus \calQ (\bd'')}$.
\end{corollary}

\begin{proof}
Immediate from Lemma~\ref{lemma intersectionbis} and
Corollary~\ref{corollary intersection}.
\end{proof}

\begin{proposition} \label{proposition diff}
Let $\lambda_0, \ldots, \lambda_p \in \bbX$ be pairwise different
and $J_l \subseteq \calA_{\lambda_l} (\bd)$, $l \in [0, p]$. If
$\bd' \in \bP + \bR$, $\bd'' \in \bQ$ and $\ol{(\calP \cup \calR)
(\bd') \oplus \calQ (\bd'')}$ is an irreducible component of
$\bigcap_{l \in [0, p]} \calH^{V_{\lambda_l, J_l}} (\bd)$, then
$\langle \bdim V_{\lambda_l, J_l}, \bd'' \rangle_{\bDelta} > 0$
for each $l \in [0, p]$. Moreover, if $\langle \bdim V_{\lambda_l,
J_l}, \bd'' \rangle_{\bDelta} = 1$ for each $l \in [0, p]$, then
there exists $M \in (\calP \cup \calR) (\bd') \oplus \calQ
(\bd'')$ such that $\partial c_{\lambda_0, J_0} (M)$, \ldots,
$\partial c_{\lambda_p, J_p} (M)$ are linearly independent.
\end{proposition}

\begin{proof}
We know from Lemma~\ref{lemma intersectionbis} that $\ol{(\calP
\cup \calR) (\bd') \oplus \calQ (\bd'')}$ is an irreducible
component of $\bigcap_{l \in [0, p]} \calH^{V_{\lambda_l}} (\bd)$.
Fix $M \in (\calP \cup \calR) (\bd') \oplus \calQ (\bd'')$ such
that $\calO (M)$ is maximal in $\bigcap_{l \in [0, p]}
\calH^{V_{\lambda_l}} (\bd)$. Write $M = P \oplus Q$ for $P \in
(\calP \cup \calR) (\bd')$ and $Q \in \calQ (\bd'')$.

First we prove that $\Hom_{\bDelta} (V_{\lambda_l, J_l}, P) = 0$
for each $l \in [0, p]$. This will imply in particular that
\[
\langle \bdim V_{\lambda_l, J_l}, \bd'' \rangle_{\bDelta} = \dim_k
\Hom_{\bDelta} (V_{\lambda_l, J_l}, Q) = \dim_k \Hom_{\bDelta}
(V_{\lambda_l, J_l}, M) > 0
\]
for each $l \in [0, p]$. Write $P = P' \oplus R$ for $P' \in \add
\calP$ and $R \in \add R$, and assume $\Hom_{\bDelta}
(V_{\lambda_l, i}, R) \neq 0$ for some $l \in [0, p]$ and $i \in
J_l$. Then $q_{\lambda_l, i}^R > 0$. If $p^R > 0$, then $\langle
\bdim Q, \bdim R \rangle_{\bDelta} \leq \langle \bd'', \bh
\rangle_{\bDelta} < 0$ by Proposition~\ref{proposition
Tits}\eqref{point Tits3} (recall that $\bd'' \neq 0$ by
Corollary~\ref{corollary intersectionbis}). Otherwise, we fix $n
\in \bbN$ such that $q_{\lambda, (i + n) \mod r_\lambda}^R = 0$
and $q_{\lambda, (i + j) \mod r_\lambda}^R > 0$ for each $j \in
[1, n - 1]$. Then
\begin{align*}
\langle \bd'', \be_{\lambda_l, i + n - 1}^n \rangle_{\bDelta} & =
\langle \bd - \bdim P' - \bdim R, \be_{\lambda_l, i + n - 1}^n
\rangle_{\bDelta}
\\
& \leq - p_{\lambda_l, (i + n) \mod r_\lambda}^{\bd} - \langle
\bdim P', \be_{\lambda, i + n - 1}^n \rangle_{\bDelta} -
q_{\lambda_l, i}^R < 0.
\end{align*}
This again implies that $\langle \bdim Q, \bdim R \rangle
_{\bDelta} < 0$, hence $\Ext_{\bDelta}^1 (Q, R) \neq 0$. If $0 \to
R \to Q' \to Q \to 0$ is a non-split exact sequence, then $P'
\oplus Q' \in \bigcap_{l \in [0, p]} \calH^{V_{\lambda_l}} (\bd)$,
since $\dim_k \Hom_{\bDelta} (V_\lambda, Q') \geq \langle \bh,
\bdim Q' \rangle_{\bDelta} = \langle \bh, \bd'' \rangle_{\bDelta}
> 0$ for each $\lambda \in \bbX$. Moreover, $M \in \ol{\calO (P'
\oplus Q')}$ and $M \not \simeq P' \oplus Q'$, which contradicts
the maximality of $\calO (M)$.

Now we assume that $\langle \bdim V_{\lambda_l, J_l}, \bd''
\rangle_{\bDelta} = 1$ for each $l \in [0, p]$ and prove that
under this assumption $\partial c_{\lambda_0, J_0} (M)$, \ldots,
$\partial c_{\lambda_p, J_p} (M)$ are linearly independent. Our
assumption implies that
\[
\dim_k \Hom_{\bDelta} (V_{\lambda_l, J_l}, M) = \dim_k
\Hom_{\bDelta} (V_{\lambda_l, J_l}, Q) = 1
\]
for each $l \in [0, p]$. Let $K := \bigcap_{l \in [0, p]} \partial
c^{V_{\lambda_l, J_l}} (M) \subseteq T_M \rep_{\bDelta} (\bd)$. We
have the canonical inclusion $\Ext_{\bDelta}^1 (Q, P)
\hookrightarrow \Ext_{\bDelta}^1 (M, M)$, which sends an exact
sequence $\xi : 0 \to P \to N \to Q \to 0$ to the sequence
\[
\xi' : 0 \to M \to N \oplus M \to M \to 0.
\]
Using Proposition~\ref{proposition Zwara}\eqref{point Zwara1} we
obtain that $\xi' \in \pi_M (K)$ if and only if $\dim_k
\Hom_{\bDelta} (V_{\lambda_l, J_l}, N) = 1$ for each $l \in [0,
p]$. In particular, this implies that $N \in \bigcap_{l \in [0,
p]} \calH^{V_{\lambda_l, J_l}} (\bd)$. By the maximality of $\calO
(M)$, $N \simeq M$, i.e.\ $\xi$ splits, thus
Proposition~\ref{proposition variety}\eqref{point variety 4}
implies
\[
\codim_{T_M \rep_{\bDelta} (\bd)} K \geq \dim_k \Ext_{\bDelta}^1
(Q, P) \geq - \langle \bd'', \bd' \rangle_{\bDelta}.
\]
It follows from Corollary~\ref{corollary inequality} that $-
\langle \bd'', \bd' \rangle_{\bDelta} \geq p + 1$, and this
finishes the proof.
\end{proof}

Let $\lambda_0, \ldots, \lambda_p \in \bbX$ be pairwise different
and $J_l \subseteq \calA_{\lambda_l} (\bd)$, $l \in [0, p]$.
Assume that $\ol{(\calP \cup \calR) (\bd') \oplus \calQ (\bd'')}$
is an irreducible component of $\bigcap_{l \in [0, p]}
\calH^{V_{\lambda_l, J_l}} (\bd)$ for $\bd' \in \bP + \bR$ and
$\bd'' \in \bQ$. We know from Corollary~\ref{corollary
intersectionbis} that $\dim (\calP \cup \calR) (\bd') \oplus \calQ
(\bd'') = a_{\bDelta} (\bd) - p - 1$. Consequently, either
$q_{\bDelta} (\bd'') = 0$ or $q_{\bDelta} (\bd'') = 1$ by
Corollary~\ref{corollary dimension}. We prove that in the latter
case there is always $M \in (\calP \cup \calR) (\bd') \oplus \calQ
(\bd'')$ such that $\partial c_{\lambda_0, J_0} (M)$, \ldots,
$\partial c_{\lambda_p, J_p} (M)$ are linearly independent.

\begin{corollary} \label{corollary independent}
Let $\lambda_0, \ldots, \lambda_p \in \bbX$ be pairwise different
and $J_l \subseteq \calA_{\lambda_l} (\bd)$, $l \in [0, p]$. If
$\bd' \in \bP + \bR$, $\bd'' \in \bQ$, $\ol{(\calP \cup \calR)
(\bd') \oplus \calQ (\bd'')}$ is an irreducible component of
$\bigcap_{l \in [0, p]} \calH^{V_{\lambda_l, J_l}} (\bd)$, and
$q_{\bDelta} (\bd'') = 1$, then there exists $M \in (\calP \cup
\calR) (\bd') \oplus \calQ (\bd'')$ such that $\partial
c_{\lambda_0, J_0} (M)$, \ldots, $\partial c_{\lambda_p, J_p} (M)$
are linearly independent.
\end{corollary}

\begin{proof}
From the previous proposition we know that $\langle \bdim
V_{\lambda_l, J_l}, \bd'' \rangle_{\bDelta} > 0$ for each $l \in
[0, p]$. On the other hand, Corollary~\ref{corollary dimension}
implies that  $\langle \bdim V_{\lambda_l, J_l}, \bd''
\rangle_{\bDelta} \leq \langle \bh, \bd'' \rangle_{\bDelta} = 1$
for each $l \in [0, p]$. Consequently, $\langle \bdim
V_{\lambda_l, J_l}, \bd'' \rangle_{\bDelta} = 1$ for each $l \in
[0, p]$, and the claim follows from the previous proposition.
\end{proof}

\section{Nonsingular dimension vectors} \label{section
nonsingular}

Throughout this section we fix a sincere separating exact
subcategory $\calR$ of $\ind \bDelta$ for a tame bound quiver
$\bDelta$ and use freely notation introduced in
Section~\ref{section tubular}. We also fix $\bd \in \bR$ such that
$p := p^{\bd} > 0$. Finally, we assume that $\bd$ is not singular.
This assumption implies, according to Proposition~\ref{proposition
singular}\eqref{point singular2} and Corollary~\ref{corollary
dimension}, that $q_{\bDelta} (\bd'') = 1$ for any $\bd' \in \bP +
\bR$ and $\bd'' \in \bQ$ such that $\bd' + \bd'' = \bd$ and $\dim
(\calP \cup \calR) (\bd') \oplus \calQ (\bd'') = a_{\bDelta} (\bd)
- p - 1$. Consequently, we have the following.

\begin{lemma} \label{lemma diffbis}
Let $\lambda_0, \ldots, \lambda_p \in \bbX$ be pairwise different
and $J_l \subseteq \calA_{\lambda_l} (\bd)$, $l \in [0, p]$. If
$\calC$ is an irreducible component of $\bigcap_{l \in [0, p]}
\calH^{V_{\lambda_l, J_l}} (\bd)$, then there exists $M \in \calC$
such that $\partial c_{\lambda_0, J_0} (M)$, \ldots, $\partial
c_{\lambda_p, J_p} (M)$ are linearly independent.
\end{lemma}

\begin{proof}
We know from Corollary~\ref{corollary intersectionbis} that $\dim
\calC = a_{\bDelta} (\bd) - p - 1$ and $\calC = \ol{\dim (\calP
\cup \calR) (\bd') \oplus \calQ (\bd'')}$ for $\bd' \in \bP + \bR$
and $\bd'' \in \bQ$. Since $q_{\bDelta} (\bd'') = 1$, the claim
follows from Corollary~\ref{corollary independent}.
\end{proof}

\begin{corollary} \label{corollary ideal}
If $\lambda_0, \ldots, \lambda_p \in \bbX$ are pairwise different,
then
\begin{multline*}
\Bigl\{ f \in k [\rep_{\bDelta} (\bd)] : \text{$f (M) = 0$ for
each $M \in \bigcap_{l \in [0, p]} \calH^{V_{\lambda_l}} (\bd)$}
\Bigr\}
\\
= (c_{\lambda_0}, \ldots, c_{\lambda_p}).
\end{multline*}
\end{corollary}

\begin{proof}
We know from Proposition~\ref{proposition variety}\eqref{point
variety 1} that $\rep_{\bDelta} (\bd)$ is a complete intersection.
Moreover, the previous lemma implies that for each irreducible
component $\calC$ of $\bigcap_{l \in [0, p]} \calH^{V_{\lambda_l}}
(\bd)$ there exists $M \in \calC$ such that $\partial
c_{\lambda_0} (M)$, \ldots, $\partial c_{\lambda_p} (M)$ are
linearly independent. Consequently, the claim follows from
Propositions~\ref{proposition Serre}\eqref{point reduced}.
\end{proof}

\begin{proposition} \label{proposition Sequiv}
Let $\lambda_0, \ldots, \lambda_p \in \bbX$ be pairwise different.
If $M, N \in \calR (\bd)$ and there exists $\mu \in k$ such that
$c_{\lambda_l} (M) = \mu c_{\lambda_l} (N)$ for each $l \in [0,
p]$, then $M$ and $N$ are S-equivalent.
\end{proposition}

\begin{proof}
Lemma~\ref{lemma intersection} implies that $c_{\lambda_l} (M)
\neq 0$ for some $l \in [0, p]$. Without loss of generality we may
assume that $c_{\lambda_0} (M) \neq 0$. Then $\mu \neq 0$ and
$c_{\lambda_0} (N) \neq 0$. For $l \in [0, p]$ we put $\mu_l :=
\frac{c_{\lambda_l} (M)}{c_{\lambda_0} (M)}$. Observe that
$c_{\lambda_l} (N) = \mu_l c_{\lambda_0} (N)$ for each $l \in [0,
p]$.

Fix $\lambda \in \bbX$, and put $\mu' := \frac{c_\lambda
(M)}{c_{\lambda_0} (M)}$ and $\mu'' := \frac{c_\lambda
(N)}{c_{\lambda_0} (N)}$. We know from Lemma~\ref{lemma
intersection} that $c_\lambda (V) = 0$ for each $V \in \bigcap_{l
\in [0, p]} \calH^{V_{\lambda_l}} (\bd)$, hence
Corollary~\ref{corollary ideal} implies that there exist $f_0,
\ldots, f_p \in k [\rep_{\bDelta} (\bd)]$ such that $c_\lambda =
\sum_{l \in [0, p]} f_l c_{\lambda_l}$. Put $f := \sum_{l \in [0,
p]} \mu_l f_l$. Then
\begin{align*}
c_\lambda (g \cdot M) & = \sum_{l \in [0, p]} f_l (g \cdot M)
c_{\lambda_l} (g \cdot M)
\\
& = \sum_{l \in [0, p]} \mu_l f_l (g \cdot M) c_{\lambda_0} (g
\cdot M) = f (g \cdot M) \cdot c_{\lambda_0} (g \cdot M)
\end{align*}
for each $g \in \GL (\bd)$. Recall that $c_\lambda$ and
$c_{\lambda_0}$ are semi-invariants of the same weight, hence $f
(g \cdot M) = \frac{c_\lambda (M)}{c_{\lambda_0} (M)} = \mu'$ for
each $g \in \GL (\bd)$. Similarly, $f (g \cdot N) = \mu''$ for
each $g \in \GL (\bd)$. Since $\ol{\calO (M)} \cap \ol{\calO (N)}
\neq \varnothing$ ($S^{\bd} \in \ol{\calO (M)} \cap \ol{\calO
(N)}$), $\mu' = \mu''$. Consequently,
\[
c_\lambda (M) = \mu' c_{\lambda_0} (M) = \mu'' \mu c_{\lambda_0}
(N) = \mu c_\lambda (N),
\]
and the claim follows from Corollary~\ref{corollary Sequivalence}.
\end{proof}

\begin{proposition} \label{proposition complete}
If $\calO (M) \subseteq \rep_{\bDelta} (\bd)$ is maximal, then
there exist $\lambda_0, \ldots, \lambda_p \in \bbX$, $i_0 \in
\calA_{\lambda_0} (\bd)$, \ldots, $i_p \in \calA_{\lambda_p}
(\bd)$, and $\mu_1, \ldots, \mu_p \in k$, such that
\begin{multline*}
\{ f \in k [\rep_{\bDelta} (\bd)] : \text{$f (N) = 0$ for each $N
\in \ol{\calO (M)}$} \}
\\
= (c_{\lambda_1, i_1} - \mu_1 c_{\lambda_0, i_0}, \ldots,
c_{\lambda_p, i_p} - \mu_p c_{\lambda_0, i_0}).
\end{multline*}
In particular, $\ol{\calO (M)}$ is a complete intersection of
dimension $a_{\bDelta} (\bd) - p$.
\end{proposition}

\begin{proof}
First, let $(\lambda_1, i_1)$, \ldots, $(\lambda_q, i_q)$ be the
pairwise different elements of $\hat{\bbX} (M)$. We put $\mu_l :=
0$ for $l \in [1, q]$. Next, we choose pairwise different
$\lambda_0, \lambda_{q + 1}, \ldots, \lambda_p \in \bbX \setminus
(\bbX_0 \cup \bbX (M))$. Finally, we put $i_0 := 0$, and $i_l :=
0$ and $\mu_l := \frac{c_{\lambda_l} (M)}{c_{\lambda_0} (M)}$ for
$l \in [q + 1, p]$.

Let
\[
\calV := \{ N \in \rep_{\bDelta} (\bd) : \text{$c_{\lambda_l, i_l}
(N) - \mu_l c_{\lambda_0, i_0} (N) = 0$ for each $l \in [1, p]$}
\}
\]
and $\calV' := \bigcap_{l \in [0, p]} \calH^{V_{\lambda_l, i_l}}
(\bd)$. Obviously $\calV' \subseteq \calV$. Moreover, every
irreducible component of $\calV'$ has dimension $a_{\bDelta} (\bd)
- p - 1$ by Corollary~\ref{corollary intersectionbis}, hence
Krull's Principal Ideal Theorem implies that every irreducible
component of $\calV$ has dimension $a_{\bDelta} (\bd) - p$. In
particular, Corollary~\ref{corollary dimension} implies that
$\calR (\bd) \cap \calC$ is a non-empty open subset of $\calC$ for
each irreducible component $\calC$ of $\calV$. Note that
$c_{\lambda_l} (R) = \frac{c_{\lambda_0} (R)}{c_{\lambda_0} (M)}
c_{\lambda_l} (M)$ for any $l \in [0, p]$ and $R \in \calR (\bd)
\cap \calV$, thus Proposition~\ref{proposition Sequiv} implies
that $R$ is S-equivalent to $M$ for each $R \in \calV \cap \calR
(\bd)$. Consequently, there are only finitely many orbits in
$\calR (\bd) \cap \calV$. This implies that every irreducible
component of $\calV$ is of the form $\ol{\calO (R)}$ for some $R
\in \calR (\bd)$. Fix $R \in \calR (\bd)$ such that $\ol{\calO
(R)}$ is an irreducible component of $\calV$. Since $\dim \calO
(R) = a_{\bDelta} (\bd) - p$, $\calO (R)$ is maximal in
$\rep_{\bDelta} (R)$. Moreover, $R$ and $M$ are S-equivalent and
$\hat{\bbX} (M) \subseteq \hat{\bbX} (R)$, hence $\calO (R) =
\calO (M)$. Consequently, $\calV = \ol{\calO (M)}$.

Lemma~\ref{lemma diffbis} implies that there exists $N \in \calV$
such that $\partial c_{\lambda_0, i_0} (N)$, \ldots, $\partial
c_{\lambda_p, i_p} (N)$ are linearly independent. Consequently,
there exists $N \in \calV$ such that $\partial c_{\lambda_1, i_1}
(N) - \mu_1 \partial c_{\lambda_0, i_0} (N)$, \ldots, $\partial
c_{\lambda_p, i_p} (N) - \mu_p \partial c_{\lambda_0, i_0} (N)$
are linearly independent. Since $\rep_{\bDelta} (\bd)$ is a
complete intersection by Proposition~\ref{proposition
variety}\eqref{point variety 1}, the claim follows from
Proposition~\ref{proposition Serre}\eqref{point reduced}.
\end{proof}

\begin{proposition} \label{proposition normal}
If $\calO (M) \subseteq \rep_{\bDelta} (\bd)$ is maximal, then the
variety $\ol{\calO (M)}$ is normal.
\end{proposition}

\begin{proof}
We know form Proposition~\ref{proposition complete} that there
exist $\lambda_0, \ldots, \lambda_p \in \bbX$, $i_0 \in
\calA_{\lambda_0} (\bd)$, \ldots, $i_p \in \calA_{\lambda_p}
(\bd)$, and $\mu_1, \ldots, \mu_p \in k$, such that
\begin{multline*}
\{ f \in k [\rep_{\bDelta} (\bd)] : \text{$f (N) = 0$ for each $N
\in \ol{\calO (M)}$} \}
\\
= (c_{\lambda_l, i_l} - \mu_l c_{\lambda_0, i_0} : \text{$l \in
[1, p]$}).
\end{multline*}
Let
\[
\calU := \{ N \in \ol{\calO (M)} : \text{$\dim_k T_N \ol{\calO
(M)} = \dim \calO (M)$} \}.
\]
Equivalently, $\calU$ is the set of all $N \in \ol{\calO (M)}$
such that $\partial c_{\lambda_1, i_1} (N) - \mu_1 \partial
c_{\lambda_0, i_0} (N)$, \ldots, $\partial c_{\lambda_p, i_p} (N)
- \mu_p
\partial c_{\lambda_0, i_0} (N)$ are linearly independent.

By general theory $\calO (M) \subseteq \calU$, hence $\ol{\calO
(M)} \setminus \calU \subseteq \calV' \cup \calV''$, where $\calV'
:= \bigcap_{l \in [0, p]} \calH^{V_{\lambda_l, i_l}} (\bd)$ and
$\calV'' := (\ol{\calO (M)} \setminus \calO (M)) \cap \calR
(\bd)$. Lemma~\ref{lemma diffbis} says that for each irreducible
component $\calC$ of $\calV'$ there exists $N \in \calC$ such that
$\partial c_{\lambda_0, i_0} (N)$, \ldots, $\partial c_{\lambda_p,
i_p} (N)$ are linearly independent. In particular, $\calU \cap
\calC \neq \varnothing$ for each each irreducible component
$\calC$ of $\calV'$, thus $\dim (\calV' \setminus \calU) < \dim
\calV' = a_{\bDelta} (\bd) - p - 1 = \dim \calO (M) - 1$. On the
other hand, if $R \in \calV''$, then $R$ is S-equivalent to $M$ by
Proposition~\ref{proposition Sequiv}, hence $\calV''$ is a union
of finitely many orbits. Moreover, \cite{Zwara2005a}*{Theorem~1.1}
implies that $R \in \calU$ for each $R \in \calV''$ such that
$\dim \calO (R) = \dim \calO (M) - 1$. Concluding, we obtain that
$\dim (\ol{\calO (M)} \setminus \calU) < \dim \calO (M) - 1$.
Since $\ol{\calO (M)}$ is a complete intersection by
Proposition~\ref{proposition complete}, the claim follows from
Proposition~\ref{proposition Serre}\eqref{point normal}.
\end{proof}

\section{Singular dimension vector} \label{section
singular}

Throughout this section we fix a sincere separating exact
subcategory $\calR$ of $\ind \bDelta$ for a tame bound quiver
$\bDelta$ and use freely notation introduced in
Section~\ref{section tubular}. We also fix singular $\bd \in \bR$.
Proposition~\ref{proposition singular}\eqref{point singular1}
implies that $\bd = \bh$ and $\bDelta$ is of type $(2, 2, 2, 2)$.
Let $\calO (M) \subseteq \rep_{\bDelta} (\bh)$ be maximal. It
follows from~\cite{BobinskiSkowronski2002}*{Proposition~5} that $M
\simeq R_{\lambda, i}^{(r_\lambda)}$ for some $\lambda \in \bbX$
and $i \in [0, r_\lambda - 1]$. We prove that $\ol{\calO (M)}$ is
normal if and only if $r_\lambda = 2$. Note that $\hat{\bbX} (M) =
\{ (\lambda, j) \}$, where $j := (i - 1) \mod r_\lambda$.
Moreover, $V_{\lambda, j} = R_{\lambda, j}$.

\begin{proposition} \label{proposition complete singular}
We have
\[
\{ f \in k [\rep_{\bDelta} (\bh)] : \text{$f (N) = 0$ for each $N
\in \ol{\calO (M)}$} \} = (c_{\lambda, j}).
\]
In particular, $\ol{\calO (M)}$ is a complete intersection of
dimension $a_{\bDelta} (\bh) - 1$.
\end{proposition}

\begin{proof}
We know from Proposition~\ref{proposition variety}\eqref{point
variety 1} that $\rep_{\bDelta} (\bh)$ is an irreducible variety
of dimension $a_{\bDelta} (\bh)$, hence Krull's Principal Ideal
Theorem implies that every irreducible component of
$\calH^{V_{\lambda, j}} (\bh)$ has dimension $a_{\bDelta} (\bh) -
1$. Observe that $\calR (\bh) \cap \calH^{V_{\lambda, j}} (\bh)$
is a union of finitely many orbits. Since $\dim
(\calH^{V_{\lambda, j}} (\bh) \setminus \calR (\bh)) \leq
a_{\bDelta} (\bh) - 2$ by Corollary~\ref{corollary dimension},
this implies that every irreducible component of $\calV$ is of the
form $\ol{\calO (R)}$ for a maximal orbit $\calO (R)$ in
$\rep_{\bDelta} (\bh)$. However,
\cite{BobinskiSkowronski2002}*{Proposition~5} implies that $\calO
(M)$ is a unique maximal orbit in $\rep_{\bDelta} (\bh)$ which is
contained in $\calH^{V_{\lambda, j}} (\bh)$, hence
$\calH^{V_{\lambda, j}} (\bh) = \ol{\calO (M)}$.

We know that $\dim_k \Ext_{\bDelta}^1 (M, M) = 1$ and the
non-split exact sequences in $\Ext_{\bDelta}^1 (M, M)$ are of the
form $\xi : 0 \to M \to N \to M \to 0$ with $N \simeq R_{\lambda,
i}^{(2 r_\lambda)}$. In particular, $\dim_k \Hom_{\bDelta}
(V_{\lambda, j}, N) = 1$. Consequently, the sequence $0 \to M \to
M \oplus M \to M \to 0$ is the only $V_{\lambda, j}$-exact
sequence in $\Ext_{\bDelta}^1 (M, M)$.
Propositions~\ref{proposition Zwara}\eqref{point Zwara1}
and~\ref{proposition variety}\eqref{point variety 4} imply that
$\partial c^{V_{\lambda, j}} (M)$ is non-zero. Since
$\rep_{\bDelta} (\bh)$ is a complete intersection by
Proposition~\ref{proposition variety}\eqref{point variety 1}, the
claim follows from Proposition~\ref{proposition Serre}\eqref{point
reduced}.
\end{proof}

\begin{proposition}
Let
\[
\calU := \{ N \in \ol{\calO (M)} : \text{$\dim_k T_N \ol{\calO
(M)} = \dim \calO (M)$} \}.
\]
\begin{enumerate}

\item
If $r_\lambda = 1$, then $\dim \ol{\calO (M)} \setminus \calU =
\dim \calO (M) - 1$. In particular, $\ol{\calO (M)}$ is not
normal.

\item
If $r_\lambda = 2$, then $\dim \ol{\calO (M)} \setminus \calU <
\dim \calO (M) - 1$. In particular, $\ol{\calO (M)}$ is normal.

\end{enumerate}
\end{proposition}

\begin{proof}
Fix $\lambda' \in \bbX \setminus (\bbX_0 \cup \{ \lambda \})$.
Lemma~\ref{lemma intersection} implies that $\rep_{\bDelta} (\bh)
\setminus \calR (\bh) = \calH^{V_\lambda} (\bh) \cap
\calH^{V_{\lambda'}} (\bh)$. By general theory $\calO (M)
\subseteq \calU$, hence $\ol{\calO (M)} \setminus \calU \subseteq
\calV' \cup \calV''$, where $\calV' := (\ol{\calO (M)} \setminus
\calO (M)) \cap \calR (\bh)$ and $\calV'' := \calH^{V_{\lambda,
j}} (\bh) \cap \calH^{V_{\lambda'}} (\bh)$. We know that $\calV'$
is a union of finitely many orbits. Moreover,
\cite{Zwara2005a}*{Theorem~1.1} implies that $R \in \calU$ for
each $R \in \calV'$ such that $\dim \calO (R) = \dim \calO (M) -
1$. Consequently, $\dim (\calV' \setminus \calU) < \dim \calV'
\leq \dim \calO (M) - 1$.

Now let $\calC$ be an irreducible component of $\calV''$.
Corollary~\ref{corollary intersectionbis} implies that $\dim \calC
= a_{\bDelta} (\bh) - 2$ and there exist $\bd' \in \bP$ and $\bd''
\in \bQ$ such that $\calC = \ol{(\calP \cup \calR) (\bd') \oplus
\calQ (\bd'')}$. Moreover, Corollary~\ref{corollary dimension}
implies that either $q_{\bDelta} (\bd'') = 1$ or $q_{\bDelta}
(\bd'') = 0$. If $q_{\bDelta} (\bd'') = 1$, then
Corollary~\ref{corollary independent} implies that $\calU \cap
\calC \neq \varnothing$.

Assume that $q_{\bDelta} (\bd'') = 0$ (according to
Proposition~\ref{proposition singular}\eqref{point singular2} this
case appears since $\bd''$ is singular). Then $\langle \bh, \bd''
\rangle_{\bDelta} = 2$ by Corollary~\ref{corollary dimension}. If
$r_\lambda = 2$, then $\langle \bdim V_{\lambda, j}, \bd''
\rangle_{\bDelta} = 1$. Indeed, we know from
Proposition~\ref{proposition diff} that $\langle \bdim V_{\lambda,
j}, \bd'' \rangle_{\bDelta} > 0$. On the other hand,
Proposition~\ref{proposition Tits}\eqref{point Tits4} implies that
$\langle \bdim V_{\lambda, j}, \bd'' \rangle_{\bDelta} = \langle
\bh, \bd'' \rangle_{\bDelta} - \langle \be_{\lambda, i}, \bd''
\rangle_{\bDelta} \leq 2 - 1 = 1$. Consequently,
Proposition~\ref{proposition diff} implies that $\calU \cap \calC
\neq \varnothing$ in this case. On the other hand, if $r_\lambda =
1$, then $\dim_k \Hom_{\bDelta} (V_{\lambda, j}, N) \geq \langle
\bh, \bd'' \rangle_{\bDelta} = 2$ for each $N \in \calC$. Thus, in
this case $\calU \cap \calC = \varnothing$ by
Proposition~\ref{proposition Zwara}\eqref{point Zwara2}.

Concluding, $\dim (\ol{\calO (M)} \setminus \calU) < \dim \calO
(M) - 1$ if and only if $r_\lambda = 2$. Since we know from
Proposition~\ref{proposition complete singular} that $\ol{\calO
(M)}$ is a complete intersection, the claims about (non-)normality
of $\ol{\calO (M)}$ follow immediately from
Proposition~\ref{proposition Serre}\eqref{point normal}.
\end{proof}

We finish this section with a remark about relationship between
the degenerations and the hom-order. Let $\bDelta'$ be a bound
quiver and $\bd_0$ a dimension vector. If $U, V \in
\rep_{\bDelta'} (\bd_0)$, then we say that $V$ is a degeneration
of $U$ (and write $U \leq_{\deg} V$) if $\calO (V) \subseteq
\ol{\calO (U)}$. Similarly, we write $U \leq_{\hom} V$ if $\dim_k
\Hom_{\bDelta'} (X, U) \leq \dim_k \Hom_{\bDelta'} (X, V)$ for
each $X \in \rep \bDelta'$ (equivalently, $\dim_k \Hom_{\bDelta'}
(U, X) \leq \dim_k \Hom_{\bDelta'} (V, X)$ for each $X \in \rep
\bDelta'$). Both $\leq_{\deg}$ and $\leq_{\hom}$ induce partial
orders in the set of the isomorphism classes of the
representations of $\bDelta'$. It is also known that $\leq_{\deg}$
implies $\leq_{\hom}$. The reverse implication is not true in
general, however $\leq_{\hom}$ implies $\leq_{\deg}$ if either
$\bDelta'$ is of finite representation type~\cite{Zwara1999} or
$\gldim \bDelta' = 1$ and $\bDelta'$ is of tame representation
type~\cite{Bongartz1995} (i.e.\ $R = \varnothing$ and $\Delta'$ is
an Euclidean quiver). We present an example showing that
$\leq_{\hom}$ does not imply $\leq_{\deg}$ for the tame concealed
canonical algebras in general.

We return to the setup of this section and assume that $r_\lambda
= 2$. Let $R := R_{\lambda, 0} \oplus R_{\lambda, 1}$. Moreover,
we fix $\bd'' \in \bQ$ such that $q_{\bDelta} (\bd'') = 0$,
$\langle \bh, \bd'' \rangle_{\bDelta} = 2$ and $\bd' \in \bP$,
where $\bd' := \bh - \bd''$. If $N \in \calP (\bd') \oplus \calQ
(\bd'')$, then
\[
\dim_k \Hom_{\bDelta} (R_{\lambda', i'}, M) \leq 1 \leq \dim_k
\Hom_{\bDelta} (R_{\lambda', i'}, N)
\]
for any $\lambda' \in \bbX$ and $i' \in [0, r_{\lambda'} - 1]$. By
adapting~\cite{BongartzFrankWolters2011}*{Corollary~4.2} to the
considered situation, we get that $R \leq_{\hom} N$ for each $N
\in \calP (\bd') \oplus \calQ (\bd'')$. On the other hand,
\[
\dim \calO (R) = a_{\bDelta} (\bd) - 2 = \dim \calP (\bd') \oplus
\calQ (\bd''),
\]
hence $\dim \calP (\bd') \oplus \calQ (\bd'') \not \subseteq
\ol{\calO (R)}$, i.e.\ there exists $N \in \dim \calP (\bd')
\oplus \calQ (\bd'')$ such that $R \not \leq_{\deg} N$.

\section{Proof of Theorem~\ref{main theorem prim}} \label{section proof}

Let $M$ be a  periodic representation of a tame
concealed-canonical quiver $\bDelta$ such that $\calO (M)$ is
maximal.

If $\Ext_{\bDelta}^1 (M, M) = 0$, then $\ol{\calO (M)} =
\rep_{\bDelta} (\bd)$ by Proposition~\ref{proposition
variety}\eqref{point variety 3}. Consequently, $\ol{\calO (M)}$ is
a normal complete intersection by Proposition~\ref{proposition
variety}\eqref{point variety 1}. Observe, that $\bdim M$ is not
singular in this case.

Now assume $\Ext_{\bDelta}^1 (M, M) \neq 0$. Using
Proposition~\ref{proposition variety}\eqref{point variety 2} we
may assume that $M \in \add \calR$ for a sincere separating exact
subcategory $\calR$ of $\ind \bDelta$.
Proposition~\ref{proposition variety}\eqref{point variety 4}
implies that $\ol{\calO (M)} \neq \rep_{\bDelta} (\bd)$.
Consequently, $p^M \neq 0$ (since $\dim \calO (M) = \dim
\rep_{\bDelta} (\bd) - p^M$) and the claim follows from
Propositions~\ref{proposition complete} and~\ref{proposition
normal}.

\bibsection


\begin{biblist}

\bib{AbeasisDelFraKraft1981}{article}{
   author={Abeasis, S.},
   author={Del Fra, A.},
   author={Kraft, H.},
   title={The geometry of representations of $A_{m}$},
   journal={Math. Ann.},
   volume={256},
   date={1981},
   number={3},
   pages={401--418},
}

\bib{AssemSimsonSkowronski2006}{book}{
   author={Assem, I.},
   author={Simson, D.},
   author={Skowro{\'n}ski, A.},
   title={Elements of the Representation Theory of Associative Algebras. Vol. 1},
   series={London Math. Soc. Stud. Texts},
   volume={65},
   publisher={Cambridge Univ. Press },
   place={Cambridge},
   date={2006},
   pages={x+458},
}

\bib{BenderBongartz2003}{article}{
   author={Bender, J.},
   author={Bongartz, K.},
   title={Minimal singularities in orbit closures of matrix pencils},
   journal={Linear Algebra Appl.},
   volume={365},
   date={2003},
   pages={13--24},
}

\bib{Bobinski2007}{article}{
   author={Bobi{\'n}ski, G.},
   title={Geometry of regular modules over canonical algebras},
   journal={Trans. Amer. Math. Soc.},
   volume={360},
   date={2008},
   number={2},
   pages={717--742},
}

\bib{BobinskiRiedtmannSkowronski2008}{article}{
   author={Bobi{\'n}ski, G.},
   author={Riedtmann, Ch.},
   author={Skowro{\'n}ski, A.},
   title={Semi-invariants of quivers and their zero sets},
   book={
      title={Trends in Representation Rheory of Algebras and Related Topics},
      editor={Skowro{\'n}ski, A.},
      series={EMS Ser. Congr. Rep.},
      publisher={Eur. Math. Soc.},
      place={Z\"urich},
   },
   date={2008},
   pages={49--99},
}

\bib{BobinskiSkowronski1999}{article}{
   author={Bobi{\'n}ski, G.},
   author={Skowro{\'n}ski, A.},
   title={Geometry of modules over tame quasi-tilted algebras},
   journal={Colloq. Math.},
   volume={79},
   date={1999},
   number={1},
   pages={85--118},
}

\bib{BobinskiSkowronski2002}{article}{
   author={Bobi{\'n}ski, G.},
   author={Skowro{\'n}ski, A.},
   title={Geometry of periodic modules over tame concealed and tubular algebras},
   journal={Algebr. Represent. Theory},
   volume={5},
   date={2002},
   number={2},
   pages={187--200},
}

\bib{BobinskiZwara2002}{article}{
   author={Bobi{\'n}ski, G.},
   author={Zwara, G.},
   title={Schubert varieties and representations of Dynkin quivers},
   journal={Colloq. Math.},
   volume={94},
   date={2002},
   number={2},
   pages={285--309},
}

\bib{BobinskiZwara2006}{article}{
   author={Bobi{\'n}ski, G.},
   author={Zwara, G.},
   title={Normality of orbit closures for directing modules over tame
   algebras},
   journal={J. Algebra},
   volume={298},
   date={2006},
   number={1},
   pages={120--133},
}

\bib{Bongartz1983}{article}{
   author={Bongartz, K.},
   title={Algebras and quadratic forms},
   journal={J. London Math. Soc. (2)},
   volume={28},
   date={1983},
   number={3},
   pages={461--469},
}

\bib{Bongartz1991}{article}{
   author={Bongartz, K.},
   title={A geometric version of the Morita equivalence},
   journal={J. Algebra},
   volume={139},
   date={1991},
   number={1},
   pages={159--171},
}

\bib{Bongartz1994}{article}{
   author={Bongartz, K.},
   title={Minimal singularities for representations of Dynkin quivers},
   journal={Comment. Math. Helv.},
   volume={69},
   date={1994},
   number={4},
   pages={575--611},
}

\bib{Bongartz1995}{article}{
   author={Bongartz, K.},
   title={Degenerations for representations of tame quivers},
   journal={Ann. Sci. \'Ecole Norm. Sup. (4)},
   volume={28},
   date={1995},
   number={5},
   pages={647--668},
   issn={0012-9593},
}

\bib{BongartzFrankWolters2011}{article}{
   author={Bongartz, K.},
   author={Frank, G.},
   author={Wolters, I.},
   title={On minimal disjoint degenerations of modules over tame path
   algebras},
   journal={Adv. Math.},
   volume={226},
   date={2011},
   number={2},
   pages={1875--1910},
   issn={0001-8708},
}

\bib{Chindris2007}{article}{
   author={Chindris, C.},
   title={On orbit closures for infinite type quivers},
   eprint={arXiv:0709.3613},
}

\bib{Chindris2009}{article}{
   author={Chindris, C.},
   title={Orbit semigroups and the representation type of quivers},
   journal={J. Pure Appl. Algebra},
   volume={213},
   date={2009},
   number={7},
   pages={1418--1429},
   issn={0022-4049},
}

\bib{CrawleyBoevey1988}{article}{
   author={Crawley-Boevey, W. W.},
   title={On tame algebras and bocses},
   journal={Proc. London Math. Soc. (3)},
   volume={56},
   date={1988},
   number={3},
   pages={451--483},
}

\bib{DerksenWeyman2002}{article}{
   author={Derksen, H.},
   author={Weyman, J.},
   title={Semi-invariants for quivers with relations},
   journal={J. Algebra},
   volume={258},
   date={2002},
   number={1},
   pages={216--227},
}

\bib{Domokos2002}{article}{
   author={Domokos, M.},
   title={Relative invariants for representations of finite dimensional algebras},
   journal={Manuscripta Math.},
   volume={108},
   date={2002},
   number={1},
   pages={123--133},
}

\bib{DomokosLenzing2002}{article}{
   author={Domokos, M.},
   author={Lenzing, H.},
   title={Moduli spaces for representations of concealed-canonical algebras},
   journal={J. Algebra},
   volume={251},
   date={2002},
   number={1},
   pages={371--394},
}

\bib{Drozd1980}{article}{
   author={Drozd, Ju. A.},
   title={Tame and wild matrix problems},
   book={
      title={Representation Theory. II},
      editor={Dlab, V.},
      editor={Gabriel, P.},
      series={Lecture Notes in Math.},
      volume={832},
      publisher={Springer},
      place={Berlin},
   },
   date={1980},
   pages={242--258},
}

\bib{Eisenbud1995}{book}{
   author={Eisenbud, D.},
   title={Commutative Algebra},
   series={Grad. Texts in Math.},
   volume={150},
   publisher={Springer},
   place={New York},
   date={1995},
   pages={xvi+785},
}

\bib{Gabriel1972}{article}{
   author={Gabriel, P.},
   title={Unzerlegbare Darstellungen. I},
   journal={Manuscripta Math.},
   volume={6},
   date={1972},
   number={1},
   pages={71--103},
}

\bib{King1994}{article}{
   author={King, A. D.},
   title={Moduli of representations of finite-dimensional algebras},
   journal={Quart. J. Math. Oxford Ser. (2)},
   volume={45},
   date={1994},
   number={180},
   pages={515--530},
}

\bib{Kunz1985}{book}{
   author={Kunz, E.},
   title={Introduction to Commutative Algebra and Algebraic Geometry},
   publisher={Birkh\"auser},
   place={Boston, MA},
   date={1985},
   pages={xi+238},
}

\bib{LenzingMeltzer1996}{article}{
   author={Lenzing, H.},
   author={Meltzer, H.},
   title={Tilting sheaves and concealed-canonical algebras},
   book={
      title={Representation Theory of Algebras},
      editor={Bautista, R.},
      editor={Mart{\'{\i}}nez-Villa, R.},
      editor={de la Pe{\~n}a, J. A.},
      series={CMS Conf. Proc.},
      volume={18},
      publisher={Amer. Math. Soc.},
      place={Providence, RI},
   },
   date={1996},
   pages={455--473},
}

\bib{LenzingdelaPena1999}{article}{
   author={Lenzing, H.},
   author={de la Pe{\~n}a, J. A.},
   title={Concealed-canonical algebras and separating tubular families},
   journal={Proc. London Math. Soc. (3)},
   volume={78},
   date={1999},
   number={3},
   pages={513--540},
}

\bib{LocZwara}{article}{
   author={Loc, N. Q.},
   author={Zwara, G.},
   title={Regular orbit closures in module varieties},
   journal={Osaka J. Math.},
   volume={44},
   date={2007},
   number={4},
   pages={945--954},
}

\bib{RiedtmannZwara2011}{article}{
   author={Riedtmann, Ch.},
   author={Zwara, G.},
   title={Orbit closures and rank schemes},
   journal={Comment. Math. Helv.},
   status={to appear},
}

\bib{Ringel1980}{article}{
   author={Ringel, C. M.},
   title={The rational invariants of the tame quivers},
   journal={Invent. Math.},
   volume={58},
   date={1980},
   number={3},
   pages={217--239},
}

\bib{Ringel1984}{book}{
   author={Ringel, C. M.},
   title={Tame Algebras and Integral Quadratic Forms},
   series={Lecture Notes in Math.},
   volume={1099},
   publisher={Springer},
   place={Berlin},
   date={1984},
   pages={xiii+376},
}

\bib{SimsonSkowronski2007a}{book}{
   author={Simson, D.},
   author={Skowro{\'n}ski, A.},
   title={Elements of the Representation Theory of Associative Algebras.
   Vol. 2},
   series={London Math. Soc. Stud. Texts},
   volume={71},
   publisher={Cambridge Univ. Press },
   place={Cambridge},
   date={2007},
   pages={xii+308},
}

\bib{SimsonSkowronski2007b}{book}{
   author={Simson, D.},
   author={Skowro{\'n}ski, A.},
   title={Elements of the Representation Theory of Associative Algebras.
   Vol. 3},
   series={London Math. Soc. Stud. Texts},
   volume={71},
   publisher={Cambridge Univ. Press },
   place={Cambridge},
   date={2007},
   pages={xii+456},
}

\bib{Skowronski1996}{article}{
   author={Skowro{\'n}ski, A.},
   title={On omnipresent tubular families of modules},
   book={
      title={Representation Theory of Algebras},
      editor={Bautista, R.},
      editor={Mart{\'{\i}}nez-Villa, R.},
      editor={de la Pe{\~n}a, J. A.},
      series={CMS Conf. Proc.},
      volume={18},
      publisher={Amer. Math. Soc.},
      place={Providence, RI},
   },
   date={1996},
   pages={641--657},
}

\bib{SkowronskiWeyman2000}{article}{
   author={Skowro{\'n}ski, A.},
   author={Weyman, J.},
   title={The algebras of semi-invariants of quivers},
   journal={Transform. Groups},
   volume={5},
   date={2000},
   number={4},
   pages={361--402},
   issn={1083-4362},
}

\bib{SkowronskiZwara1998}{article}{
   author={Skowro{\'n}ski, A.},
   author={Zwara, G.},
   title={Degenerations for indecomposable modules and tame algebras},
   journal={Ann. Sci. \'Ecole Norm. Sup. (4)},
   volume={31},
   date={1998},
   number={2},
   pages={153--180},
   issn={0012-9593},
}

\bib{SkowronskiZwara2003}{article}{
   author={Skowro{\'n}ski, A.},
   author={Zwara, G.},
   title={Derived equivalences of selfinjective algebras preserve
   singularities},
   journal={Manuscripta Math.},
   volume={112},
   date={2003},
   number={2},
   pages={221--230},
}

\bib{Zwara1999}{article}{
   author={Zwara, G.},
   title={Degenerations for modules over representation-finite algebras},
   journal={Proc. Amer. Math. Soc.},
   volume={127},
   date={1999},
   number={5},
   pages={1313--1322},
   issn={0002-9939},
}

\bib{Zwara2002}{article}{
   author={Zwara, G.},
   title={Unibranch orbit closures in module varieties},
   journal={Ann. Sci. \'Ecole Norm. Sup. (4)},
   volume={35},
   date={2002},
   number={6},
   pages={877--895},
}

\bib{Zwara2003}{article}{
   author={Zwara, G.},
   title={An orbit closure for a representation of the Kronecker quiver with
   bad singularities},
   journal={Colloq. Math.},
   volume={97},
   date={2003},
   number={1},
   pages={81--86},
}

\bib{Zwara2005a}{article}{
   author={Zwara, G.},
   title={Regularity in codimension one of orbit closures in module
   varieties},
   journal={J. Algebra},
   volume={283},
   date={2005},
   number={2},
   pages={821--848},
}

\bib{Zwara2005b}{article}{
   author={Zwara, G.},
   title={Orbit closures for representations of Dynkin quivers are regular
   in codimension two},
   journal={J. Math. Soc. Japan},
   volume={57},
   date={2005},
   number={3},
   pages={859--880},
}

\bib{Zwara2006}{article}{
   author={Zwara, G.},
   title={Singularities of orbit closures in module varieties and cones over
   rational normal curves},
   journal={J. London Math. Soc. (2)},
   volume={74},
   date={2006},
   number={3},
   pages={623--638},
}

\bib{Zwara2007}{article}{
   author={Zwara, G.},
   title={Codimension two singularities for representations of extended
   Dynkin quivers},
   journal={Manuscripta Math.},
   volume={123},
   date={2007},
   number={3},
   pages={237--249},
}

\bib{Zwara2011}{article}{
   author={Zwara, G.},
   title={Singularities of orbit closures in module varieties},
   book={
      title={Representations of Algebras and Related Topics},
      editor={Skowro{\'n}ski, A.},
      editor={Yamagata, K.},
      series={EMS Ser. Congr. Rep.},
      publisher={Eur. Math. Soc.},
      place={Z\"urich},
   },
   date={2011},

}

\end{biblist}
\end{document}